\documentclass[reqno, final]{amsart}
\usepackage{amsmath, amsfonts}
\usepackage{hyperref}

\newtheorem{lemma}{Lemma}[section]

\newtheorem{theorem}[lemma]{Theorem}

\newtheorem{remark}[lemma]{Remark}
\newtheorem{proposition}[lemma]{Proposition}

\newcommand{\R}{{\mathbb R}}
\newcommand{\N}{{\mathbb N}}
\newcommand{\E}{{\mathbb E}}
\newcommand{\PP}{{\mathbb P}}
\newcommand{\T}{\theta}
\numberwithin{equation}{section}

\begin{document}

\title[Large deviations for the Hermitian Brownian motion]{Large deviations for the largest eigenvalue of an Hermitian Brownian motion}

\author{C. Donati-Martin}
\author{M. Ma\"\i da}

\address{Universit\'e Versailles-Saint Quentin\newline
  Laboratoire de Math\'ematiques \newline
45 av. des Etats Unis \newline 
 F-78035 Versailles cedex, France.}
\email{catherine.donati-martin@uvsq.fr} 
\urladdr{\url{http://lmv.math.cnrs.fr/annuaire/donati-martin-catherine/}}

\address{Universit\'e Paris Sud \newline
Laboratoire de Math\'ematiques \newline
 Facult\'e des Sciences \newline
 91405 Orsay Cedex, France.}
\email{mylene.maida@math.u-psud.fr}
\urladdr{\url{http://www.math.u-psud.fr/~maida/}}

\thanks{Research supported by  \emph{Agence Nationale de la Recherche} grants ANR-08-BLAN-0311-03 and ANR-09-BLAN-0084-01}


\begin{abstract}
We establish a large deviation principle for  the process of the largest eigenvalue of an Hermitian Brownian motion.  
By a contraction principle, we recover the LDP for the largest eigenvalue of a rank one deformation of the GUE.
\end{abstract}

\maketitle

%%%%%%%%%%%%%%%%%%%%%%%%%%%%%%%%%%%%%%%%%%%%%%%
%% Author commands and definitions
%%%%%%%%%%%%%%%%%%%%%%%%%%%%%%%%%%%%%%%%%%%%%%%

%%%%%%%%%%%%%%%%%%%%%%%%%%%%%%%%%%%%%%%%%%%%%%%%

%\date{March 14, 2011; accepted September 21, 2012}

\section{Introduction}

The Gaussian unitary ensemble (GUE) is probably the most studied ensemble of random matrices.
In this work, we will focus on a dynamical version of the GUE, introduced in 1962 by Dyson : he defined
the Hermitian Brownian motion whose set of eigenvalues is a time-dependent Coulomb gas,
 consisting in particles evolving according to Brownian motions under the influence of their mutual electrostatic repulsions. 
More precisely, 
% Let us recall the definition of Dyson's Brownian motion.
%Let  $(H_N^\T(t))_{0\leq t \leq 1}$ a   Hermitian Brownian motion defined on a probability space $(\Omega, {\mathcal F}, \PP)$, 
%starting from  $H_N^\T(0) := {\rm{diag}}(\theta, 0, \ldots, 0),$ with $\T \geq 0.$
let $(\beta_{ij}, \beta^\prime_{ij})_{1 \le i \le j \le N}$ be a collection of independent identically distributed  standard real 
Brownian motions; %defined on a probability space $(\Omega, {\mathcal F}, \PP)$; 
 the Hermitian Brownian motion $(H_N(t))_{t \geq 0}$ is the random process, taking values in the 
space of $N \times N$ Hermitian matrices, with entries $(H_N)_{kl}$ given, for $k \le l$ by
\[
(H_N)_{kl} = \left\{ \begin{array}{ll}
                      \frac{1}{\sqrt{2N}} (\beta_{kl} + i \beta^\prime_{kl}), & \textrm{ if } k <l,  \\
 \frac{1}{\sqrt{N}} \beta_{kk}, & \textrm{ if } k=l.
                     \end{array}
\right.
\]
\cite{D} (see also \cite[12.1]{G}) showed that the eigenvalues of $(H_N(t))_{t \ge 0}$ satisfy the following system of stochastic differential equations (SDE) 
\begin{equation} \label{EDSvp}
d\lambda_i(t) = \frac{1}{\sqrt{N}}dB_i(t) + \frac{1}{N} \sum_{ j \not= i} \frac{1}{\lambda_i(t) - \lambda_j(t)} dt , \, t \geq 0, \; i,j = 1, \ldots, N
\end{equation}
where $B_i$ are independent standard real  Brownian motions.

It was rigorously shown in \cite{CL} that this system of SDE  admits a unique strong solution  and the
eigenvalues do not collide.

The process of the eigenvalues is called Dyson Brownian motion. Almost surely (a.s.), for any $t \ge 0,$ the corresponding  spectral measure 
$(\mu_N)_t :=\frac{1}{N} \sum_{i = 1}^N \delta_{\lambda_i(t)}$  converges weakly to the
  semicircular distribution $\sigma_t$ given by 
\begin{equation} \label{semicirc}
 d\sigma_t(x) =\frac{1}{2 \pi t} {\bf 1}_{[-2\sqrt t, 2\sqrt t]} \sqrt{4t-x^2} dx
\quad \textrm{ and } \quad \sigma_0 = \delta_0.
\end{equation}
Let us now recall some large deviations results. For the global regime in the static case (GUE), a large deviation principle (LDP) in the scale $N^2$ was  established in \cite{BAG}  
for the spectral measure of size $N.$ There exists also a dynamical version of this result :  
\cite{CDG}, \cite{G} showed the following.

Let $C([0,1]; {\mathcal
    P}(\mathbb{R})),$ the set of continuous functions on $[0,1]$ with values in the set ${\mathcal
    P}(\mathbb{R})$ of probability measures on $\R$. We equip this set with the metric
$d(\mu, \nu) = \sup_{t \in [0,1]} d_{Lip} (\mu_t, \nu_t)$ where
 $$d_{Lip}(\mu_t, \nu_t) = \sup_{f \in {\mathcal F}_{Lip}} \left| \int f d\mu_t - \int f d\nu_t \right|$$
 where ${\mathcal F}_{Lip}$ denotes the space of bounded Lipschitz functions on $\mathbb{R}$ with Lipschitz and uniform bound less than 1. 
Then 
the process $((\mu_N)_t)_{0 \le t \le 1}$
  satisfies a LDP in the scale $N^2$  with respect to the topology inherited from the metric $d.$  \smallskip
%For any $\mu \in  C([0,1]; {\mathcal P}(\mathbb{R}))$ and $\alpha > 0,$ $\mathbb B(\mu, \alpha)$ will denote the ball centered at $\mu$ with radius
%$\alpha$ with respect to the metric $d.$ \\ 

We now define $H_N^\T(t ) = H_N(t) + H_N^\T(0)$ the Hermitian Brownian motion starting from  $H_N^\T(0) := {\rm{diag}}(\theta, 0, \ldots, 0),$ with $\T \geq 0$ and
%normalized by $\E[ |(H_N^\T(t)-H_N^\T(0))_{ij}|^2] = \frac{t}{N} , 1\leq i,j \leq N$.
 denote by $ \lambda_1^{\T,N}(t)  \geq \lambda_2^{\T,N}(t) \geq \ldots \geq
  \lambda_N^{\T,N}(t)$ the set of ordered eigenvalues of
  $H_N^\T(t)$. All the results stated above about the global regime of the Hermitian Brownian motion starting from 0 will stay 
valid for $(H_N^\T(t))_{t \geq 0}.$ In this work, we will be interested in the process of the maximal eigenvalue $(\lambda_1^{\T,N}(t))_{t \geq 0},$ that is the largest particle of Dyson Brownian motion.

In the case $\T = 0,$  the corresponding quantity in the static case is just the maximal eigenvalue of the GUE.
It is well known (see for example \cite{BY}) that it converges a.s. to 2.
In  the case $\T >0$  a similar result holds  for a rank one additive deformations of the GUE  
(see for example \cite{P}).\smallskip

These results can be easily extended to  the a.s. convergence of our process, in the topology of uniform convergence for continuous 
functions on $[0,1],$ towards the function $(f_\T(t))_{ t \geq 0}$ given by:
  $$ \left \{ \begin{array}{ll}
  f_\T(t) = 2\sqrt{t} & \mbox{if $\theta = 0$,} \\
   f_\T(t) = \left\{ \begin{array}{ll}  \T + \frac{t}{\T} & t \leq \theta^2 \\
   2\sqrt{t} & t \geq \T^2
   \end{array} \right.& \mbox{if $\theta > 0$.} 
  \end{array} \right.$$
In particular, this result can be seen as a direct consequence of our main result stated below.

At the level of large deviations, the LDP for the largest eigenvalue of the GOE with a scale $N$ was obtained 
in \cite{BDG}. This result was extended to a rank one perturbation of the GUE/GOE by one of the author in \cite{M}. 
The main goal of this paper will be to prove a dynamical version of these two results. More precisely, we will consider the process
$(\lambda_1^{\T,N}(t))_{ 0\leq t \leq 1}$ as a sequence of random variables with values in the space $C_\T([0,1], \R)$ of continuous functions 
from $[0,1]$ to $\R$ equal to $\T$ at zero and investigate its LDP in this space endowed with the uniform convergence.
Our main result can be stated as follows.

\begin{theorem} \label{main}
%Let $(\lambda_1^{\T,N}(t))_{ 0\leq t \leq 1}$ be  the process of the largest
%igenvalue of an Hermitian Brownian motion  $(H_N^\T(t))_{0\leq t \leq
  %1}.$  Then  
The law of  $(\lambda_1^{\T,N}(t))_{ 0\leq t \leq
  1}$
satisfies a large deviation principle on  $C_\T([0,1]; \mathbb{R})$ equipped with the topology of uniform convergence,
 in the scale $N,$ with good rate function 
\begin{equation} \label{rf}
 I_\T(\varphi) = \left\{ \begin{array}{ll}
\displaystyle \frac 1 2\int_0^1 \left( \dot{\varphi}(s) -\frac{1}{2s}\left(\varphi(s) -
\sqrt{\varphi(s)^2-4s}\right)\right)^2 ds, \\
 \textrm{ if } \varphi \textrm{ absolutely continuous and } \varphi(t) \geq 2\sqrt{t} \; \forall t \in [0,1],\\
+ \infty,    \textrm{ otherwise. }
\end{array}
\right.
\end{equation} 
\end{theorem}

As a consequence of our main result, we will recover by contraction the fixed-time large deviation principles already shown in \cite{BDG} and  \cite{M}. The following result is a corrected version of Theorem 1.1 in \cite{M}, the proof there is correct.
 
\begin{theorem}
\label{pgdflou}
The largest eigenvalue of $H_N^\T(1)$ 
satisfies a large deviation principle in the scale $N,$
with good rate function $K_\theta$  defined as follows:

\smallskip
\noindent$\bullet$ If $\theta\leq 1,$ 
\[
K_\theta(x) = \left\{
\begin{array}{ll}
+\infty, & \textrm{ if } x< 2\\
\displaystyle \int_{2}^x \sqrt{z^2-4} \,dz, & \textrm{ if } 2 \leq x
\leq \theta +\frac 1{\theta}, \\
M_\theta(x), &  \textrm{ if }  x \geq  \theta +\frac 1{\theta},

\end{array}
\right.\]
\[
\textrm{with }  \,\,\displaystyle  M_\theta(x) = \frac 1 2 \int_{2}^x
\sqrt{z^2-4} \,dz
-\theta x +\frac 1 4 x^2 + \frac 1 2 + \frac
1 2 \theta^2 + \log \theta. \]

\noindent $\bullet$ If $\theta\geq 1,$ 
\[
K_\theta(x) = \left\{
\begin{array}{ll}
+\infty, & \textrm{ if } x< 2\\
L_\theta(x), &  \textrm{ if }  x \geq 2,
\end{array}
\right.
\]
\[
\textrm{with }\,\, \displaystyle  L_\theta(x) = \frac 1 2 \int_{ \theta +\frac 1{\theta}}^x
\sqrt{z^2-4} dz - \theta\left(x-\left( \theta +\frac 1{\theta}\right)\right) 
+ \frac 1 4 \left(x^2-\left( \theta +\frac 1{\theta}\right)^2\right).\]
\end{theorem}
\noindent

Before going further, let us make a few remarks :

\begin{remark}
\begin{enumerate}
 \item  For the sake of simplicity the theorems above are stated and proven in the paper for the Hermitian Brownian motion but we want to mention
that they can be easily extended to the symmetric Brownian motion. With the notations already introduced above, the latter is defined as the random process taking values in the space of $N \times N$ real symmetric matrices so that
\[
(S_N)_{kl} =         \frac{1}{\sqrt{N}} \beta_{kl} ,  \textrm{ if } k < l, \, (S_N)_{kk} =       \sqrt{  \frac{2}{N}}\beta_{kk} .\]  
The process of its eigenvalues satisfies the following  system of SDE
$$
d\lambda_i(t) = \frac{\sqrt 2}{\sqrt{N}}dB_i(t) + \frac{1}{N} \sum_{ j \not= i} \frac{1}{\lambda_i(t) - \lambda_j(t)} dt , \, t \geq 0, \; i = 1, \ldots, N$$
and its law satisfies a LDP with good rate function simply given by
$\frac{1}{2}I_\T.$

\smallskip\noindent
The proof that will be developed in the sequel can be adapted to the symmetric case with a few minor changes left to the reader.
\item In the whole introduction till now, we have considered our processes on the set of times $t \in [0,1]$ for simplicity but everything
could be easily generalised to any compact set $[0,T]$ for $T>0.$
\item In the sequel, we will specify the superscript $\T,N$ in the statements but drop it in the proofs,
unless there is any ambiguity. 
\end{enumerate}
\end{remark}

Let us now specify a little bit the main features of the strategy of the proof. 
In both  \cite{BDG} and  \cite{M}, the fact that the deviations of the spectral measure and those of the largest eigenvalues do not occur in the same scale plays a crucial role
and so will be in the proof of our result : in the scale at which we look at the largest eigenvalue, the spectral measure of all but the largest eigenvalue
is already well concentrated around the semicircle law.

In the static case, the LDP was shown using the explicit expression of the joint distribution of the $N$ eigenvalues. In the dynamical case, the proof relies on stochastic calculus
using that the  process of the eigenvalues satisfies the system of SDE  \eqref{EDSvp}. 
Roughly speaking, the largest eigenvalue is a solution of a SDE of the form
$$d\lambda_1(t) = \frac{1}{\sqrt{N}}dB(t) + b(\lambda_1(t), (\nu_N)_t) dt,$$
with $B$ a standard real Brownian motion, 
$\nu_N := \frac{1}{N-1} \sum_{i = 2}^N \delta_{\lambda_i(t)}$  the empirical distribution of all but the largest eigenvalues
and  the drift $b(x, \nu) $ to be explicited in the sequel. 
In the scale of interest, $\nu_N$ is close to $\sigma$ and the rate function $I_\T$ is the one predicted by 
the Freidlin-Wentzell Theorem (see \cite[Th.5.6.3]{DZ}) for the SDE
$$d\lambda_1(t) = \frac{1}{\sqrt{N}}dB(t) + b(\lambda_1(t), \sigma_t) dt.$$
One of the main difficulties will be to deal with the singularity of the drift $b,$ as for some $x \in \R,$ $\nu \mapsto b(x,\nu)$ 
is not a continuous function for the weak convergence of probabilities.\smallskip

%\noindent
The organisation of the paper will be the following.
To prove our main result, we  first establish the exponential tightness of the process $(\lambda_1(t))_{0\le t \le 1}$ stated in Proposition \ref{exptight} and proved 
in Section \ref{sec:tight}.
% and the correct lower bound,  stated in Proposition \ref{lowerbound} and proved in Section \ref{sec:lower}.
A  short Section \ref{sec:rf} will be devoted to the study of the rate function $I_\T,$ where we check in particular its lower semicontinuity.
Section \ref{sec:lower} is devoted to the proof of the lower bound, stated in Proposition \ref{lowerbound0}. The upper bound is given in 
\eqref{wub3} and
obtained along Section \ref{sec:upper}.
Then  Theorem \ref{main} will follow from the exponential tightness, 
the lower bound obtained in Proposition \ref{lowerbound0} and the weak upper bound \eqref{wub3} (see \cite[Chapt. 4]{DZ}
or \cite[Cor. D.6 and Th. D.4]{AGZ}).
 Finally, in Section \ref{sec:contract}, we recover by contraction principle the fixed-time LDP stated in Theorem \ref{pgdflou}.
%But we could not get the correct upper bound. Nevertheless, Section \ref{sec:upper} gathers some partial results in this direction.\\
%Another clue that Conjecture \ref{main} above should be true is given in Section \ref{sec:contract} : under the assumption that it holds,
% is correct. \\

%The organization of the paper is as follows: in Section 2, we prove that the law of the largest eigenvalue process is exponentially tight. In Section 3, we prove the weak upper bound. Section 4 is devoted to the proof of the lower bound. Finally, in Section 5, by applying a contraction principle, we recover the static case ($t=1$) by computing the rate function $K_\T$  from the rate function $I_\T$. \\

%

 % 

 %Finally, for $\mu \in \mathcal P(\R),$ we denote by $r(\mu)$ the right
 % end-point of the support of $\mu.$

%%%%%%%%%%%%%%%%%%%%%%%%%%%%%%%%%%%%
\section{Exponential tightness}
 \label{sec:tight}

 We want to show the exponential tightness of the process $(\lambda_1(t))_{0\le t \leq 1}$ in scale $N$ that is
\begin{proposition}\label{exptight}
 For all L, there exists $N_0$ and a compact set $K_L$ of $C_\T([0,1]; \R)$ such that:
 $$\forall N \geq N_0, \PP(\lambda_1^{\T,N} \not\in K_L) \leq \exp(-LN).$$
 \end{proposition}

 \noindent
 From the description of the compact sets of $C([0,1]; \R)$ (Ascoli theorem), it is enough to show  (see \cite[Chapter XIII, Section 1]{RY}, \cite[Section 2.3]{CDG}) the following lemma 
\begin{lemma}
 \label{lemexptight} 
 For any  $\eta >0$, there exists $\delta_0$ such that for any $\delta < \delta_0,$ for all $N,$ $p \le N$ and $s \in [0,1],$
$$\PP\left(\sup_{ s \leq t\leq s+\delta }  | \lambda^{\T,N}_p(t) - \lambda^{\T,N}_p(s) | \geq \eta\right) \leq \exp\left(-\frac{1}{10} N \frac{\eta^2}{\delta}\right). $$
\end{lemma}

\noindent
To get the proposition, for a fixed $L,$ we choose $p=1,$ any $\eta$ and then $\delta$ small enough so that $\frac{\eta^2}{10\delta} > L.$

 %$ \forall L,\,  \forall \eta,\,  \forall s \in[0,1]$,  there exists a number $\delta$, $0< \delta<1$ such that, for $N$ large enough, 
% $$\PP\left(\sup_{ s \leq t\leq s+\delta }  | \lambda_1(t) - \lambda_1(s) | \geq \eta\right) \leq \exp(-LN). $$
 
 \noindent
{\bf Proof of lemma \ref{lemexptight}}: Let $0 \leq s \leq 1$. 

\smallskip\noindent
Let us denote by $\tilde{H}_N$ the Hermitian Brownian motion defined, for $u \ge 0,$ by
$\tilde{H}_N(u) = H_N^\T(u+s ) - H_N^\T(s)$  and by
$(\tilde{\lambda}_i (u))_{u \ge 0}$ its eigenvalues, in decreasing
order.
 From a classical relation between eigenvalues (usually called Weyl's interlacing inequalities), for $t \geq s$, 
$$\lambda_p^{\T,N} (s) + \tilde{\lambda}_N(t-s) \leq \lambda_p^{\T,N} (t) \leq \lambda_p^{\T,N} (s) + \tilde{\lambda}_1(t-s) $$
so that

 $$ \vert \lambda_p^{\T,N} (t) - \lambda_p^{\T,N} (s) \vert  \leq \max\left( \tilde{\lambda}_1(t-s), -  \tilde{\lambda}_N(t-s)\right) =  \Vert \tilde{H}_N(t-s) \Vert $$
 where $\Vert . \Vert$ denotes the operator norm on matrices.  For any  $\eta >0$

\begin{eqnarray*} 
\PP\left(\sup_{ s \leq t\leq s+\delta }  | \lambda_p^{\T,N} (t) - \lambda_p^{\T,N} (s) | \geq \eta\right) &\leq& \PP\left(\sup_{ s \leq t\leq s+\delta } \Vert \tilde{H}_N(t-s)\Vert  \geq \eta\right) \\
&=& \PP\left(\sup_{ 0 \leq u \leq \delta }  \Vert \tilde{H}_N(u) \Vert \geq \eta\right) \\
&=& \PP\left(\sqrt{\delta} \sup_{ 0 \leq u \leq 1}  \Vert H_N^0(u) \Vert \geq \eta\right), \\
%&\leq & \exp(-LN) E[\exp( NL\sqrt{\delta}/ \eta \, \sup_{0\leq u \leq 1} \Vert H_N(u) \Vert)]
\end{eqnarray*}
where we used that $\tilde{H}_N$ has the same law as $H_N^0$ and the scaling invariance of this law.
Therefore, it is enough to show that for $M$ large enough and for all $N$,
\begin{equation} \label{concentration}
 \PP \left( \sup_{ 0 \leq u \leq 1} \Vert H_N^0(u) \Vert \geq M\right) \leq \exp( - N M^2/10).
 \end{equation}
Indeed, from \cite{BDG}, Lemma 6.3, we have: for $y$ large enough, for all $N,$
\begin{equation} \label{concGUE}
\PP( \Vert H_N^0(1) \Vert \geq y) \leq \exp( -N y^2/9).
\end{equation}
This implies that for $0<a < 1/9$, $\E(e^{aN \Vert H_N^0(1) \Vert ^2}) \leq 2 e^{aNy_0^2}$ for some $y_0$ large enough  and all $N$.
Now, $( \Vert H_N^0(u) \Vert )_{ 0 \leq u\leq 1}$ is a positive submartingale and from Doob's inequalities, all the moments of $ \sup_{0 \leq u \leq1}\Vert H_N^0(u) \Vert $
 are bounded  by those of $ \Vert H_N^0(1) \Vert $ (up to a constant 4). Therefore, 
 $$ \E(e^{aN \sup_{ 0 \leq u \leq 1} \Vert H_N^0(1) \Vert ^2}) \leq 4 \E(e^{aN \Vert H_N^0(1) \Vert^2}) \leq 8 e^{aNy_0^2}.$$
 From Markov's inequality,
 $$ \E( \sup_{ 0 \leq u \leq 1} \Vert H_N^0(1) \Vert  \geq M) \leq 8 e^{-aNM^2} e^{aNy_0^2} \leq e^{-a'NM^2} $$
 for $0<a'<a$ and $M$ large enough, proving \eqref{concentration}. 
 %The Gaussian concentration inequality \eqref{concGUE} implies the same type of
% concentration for $\sup_{0 \leq u \leq 1}\Vert H_N^0(u) \Vert, $ that is \eqref{concentration}.
\hfill $\Box$

%%%%%%%%%%%%%%%%%%%%%%%%%%%%%%%%%%%%%%%%%%%%%%%%%%%%%%%%%%%%%%%%%%%%
\section{Some insight on the expected rate function}
\label{sec:rf}

Before going into the proof of the lower bound, we gather hereafter some useful remarks about the function $I_\T$
defined in Theorem \ref{main}. In particular, we show in this section that it is lower semi-continuous.

We introduce the following notations :
 for $\mu $ a 
 %compactly supported 
 probability measure on $\R$ and $x \in \R$,
 %which is greater than $r(\mu)$ the right end point of the support of $\mu$
  we define
 \begin{equation} \label{defdrift}
 b(x, \mu) = \int_{-\infty}^x  \frac{ d\mu(y)}{x-y}  
 %:= \lim_{ \varepsilon \vers 0} \int_{-\infty}^{x-\varepsilon} \frac{ d\mu(y)}{x-y}  
 \in \mathbb{R}_+ \cup \{ \infty\}.
 \end{equation}
% Note that the drift term in \eqref{EDSvp} for $\lambda_1$ can also be written as $\frac{N-1}{N} b(\lambda_1(t), \nu_N(t))$. 
 
\noindent
 For $\mu \in \mathcal P(\R),$ we denote by $r(\mu)$ the right
end-point of the support of $\mu.$
Let  $(\varphi, \mu) \in
 C_\T([0,1]; \R) \times C([0,1]; \mathcal P(\R)) $ such that for all $t
 \in [0,1],$  $\varphi(t) > r(\mu_t).$ Then, $b(\varphi(t), \mu_t)$ is bounded.  We set
$$\mathcal H := \{h  \in C([0,1], \R) / h \textrm{ absolutely continuous}, \, \dot{h} \in \mathbb L^2([0,1]) \}$$ with 
$\mathbb L^2([0,1])$ the set of square-integrable functions from $[0,1]$ to $\R$
equipped with its usual $\mathbb L^2$-norm, denoted by $\Vert \cdot\Vert_2$. For any $h \in \mathcal H$ we   define
  \begin{eqnarray}
  G(\varphi, \mu; h) &=& h(1) \varphi(1)  - h(0) \varphi(0) -
   \int_0^1 \varphi(s) \dot{h}(s) ds - \int_0^1 b(\varphi(s), \mu_s)
   h(s) ds,  \label{defG0} \nonumber\\
   F(\varphi, \mu; h)  &:=& G(\varphi, \mu; h) - \frac{1}{2} \int_0^1 h^2(s) ds. \label{defF}\
  \end{eqnarray}
%Note that $F(\varphi, \sigma; h)$ is also well defined for $\varphi(t) \geq 2 \sqrt{t}$.
 
 For $\sigma := (\sigma_t)_{t\ge 0}$ the semicircular process defined in \eqref{semicirc}, the condition  
$\varphi(t) > r(\mu_t)$ reads $\varphi(t) > 2\sqrt t$ and one can check that $F(\varphi, \sigma; h)$ is also well 
defined under the weaker assumption that $\varphi(t) \ge 2\sqrt t$ for all $t \in [0,1].$ It is indeed well known (see for example \cite[p. 94]{HP})
that  
\begin{equation}
 \label{Hilbert}
b(\varphi(t) , \sigma_t  )= \frac{1}{2t} ( \varphi(t)- \sqrt{ \varphi^2(t) - 4t}), 
\end{equation}
 so that $0 \leq b(\varphi(t) , \sigma_t) \leq \frac{1}{\sqrt{t}}$ for $\varphi(t) \geq 2 \sqrt{t}.$\smallskip

\noindent
We now study the properties of $F.$
%and     $\mathcal H_- := \{ h \in \mathcal H / h \leq 0\}.$

\begin{lemma} \label{IetJ}
Let $\T \geq 0$ and $\varphi \in  C_\T([0,1], \R) $ such that for any $t \in [0,1],$
$\varphi(t) \geq 2 \sqrt t$ and define 
\begin{equation} \label{fcttaux}
 J(\varphi) := \sup_{h\in \mathcal H} F(\varphi, \sigma; h).
 % \textrm{ and } I_-(\varphi)  = \sup_{h \in \mathcal H_-} F(\varphi, \sigma; h).
 \end{equation}
Then,
\begin{itemize}
\item[]  \begin{center} 
$J(\varphi) < \infty \Rightarrow \varphi$ absolutely
  continuous \\ and  $\displaystyle J(\varphi) = \frac{1}{2} \int_0^1
  (\dot{\varphi}(s) - b(\varphi(s), \sigma_s))^2 ds = I_\T(\varphi).$          
         \end{center}
\end{itemize}
\end{lemma}

 \noindent{\bf  Proof:}
% The optimization problem for $J$ is a classical one. We just give an outline of the proof
%and refer for example to \cite{DRYZ} for details on similar computations.
Recall that $F(\varphi, \sigma; h) = G(\varphi, \sigma; h) -
 \frac{1}{2} \int_0^1 h^2(s) ds$ where $h \mapsto G(\varphi, \sigma;
 h)$ is a linear functional (see \eqref{defG0}). Replacing $h$ by $\lambda h$, $\lambda \in \R $ and optimizing in $\lambda$ yields
$$ J(\varphi) = \frac{1}{2} \sup_{h \in \mathcal H} \frac{G^2(\varphi, \sigma; h)}{ \Vert h \Vert^2_{2}}. $$ 
If $J(\varphi) < \infty$, then the linear functional $G(\varphi,
\sigma, .)$ can be extended continuously to $\mathbb L^2([0,1])$ and by Riesz
theorem, there exists $k_\varphi \in \mathbb L^2([0,1])$ such that $G(\varphi,
\sigma,h) = \int_0^1 h(s) k_\varphi(s) ds$. Comparing with \eqref{defG0}, we see  that 
$\varphi - \int_0^. b(\varphi(s), \mu_s)ds $ is absolutely continuous and  $k_\varphi(s) =
\dot{\varphi}(s) - b(\varphi(s), \sigma_s).$ 

 From Cauchy-Schwarz inequality, we obtain:
$$ G^2(\varphi, \sigma; h) \leq || k_\varphi||_2^2 ||h||_2^2  $$
with equality if $h$ is proportional to $k_\varphi$.
Therefore, $J(\varphi) \leq  \frac{1}{2} || k_\varphi||_2^2$, and the equality holds since 
 $\mathcal H$ is dense in $\mathbb L^2([0,1])$.
 %$J(\varphi) =  \frac{1}{2} \Vert k_\varphi \Vert_{2}^2$. 
The  equality  between $\frac{1}{2} || k_\varphi||^2$ and $I_{\T}(\varphi)$  follows from the computation of the Hilbert transform 
of the semicircular distribution recalled in \eqref{Hilbert} and $ \varphi(t) \geq 2 \sqrt{t}$.
\hfill  $\Box$

\noindent
We can now show the following :
\begin{proposition} \label{sci}
The function $I_\T: C_\T([0,1], \R) \rightarrow \R$  is lower semicontinuous.
\end{proposition}

 \noindent{\bf Proof:}
 From Lemma \ref{IetJ}, $I_\T(\varphi) = \sup_{h \in  \mathcal H} F(\varphi, \sigma; h)$ where 
\begin{multline*}
 F(\varphi, \sigma; h) = h(1) \varphi(1)  - h(0) \varphi(0) - \int_0^1 \varphi(s) \dot{h}(s) ds - \int_0^1 b(\varphi(s), \sigma_s) h(s) ds 
  \\
- \frac{1}{2} \int_0^1 h^2(s) ds.
\end{multline*}
 % \end{eqnarray*}
  We shall prove that for fixed $h \in {\mathcal H}$, $\varphi \mapsto F(\varphi, \sigma; h)$ is continuous. From the definition of $F$, performing an integration
 by part in the term of the integral in $b$, it is enough to prove the continuity of $\varphi \mapsto \Lambda(\varphi) := \int_0^. b(\varphi(s), \sigma_s) ds$ 
(the other terms are obviously continuous in $\varphi$).
As we know that $0 \leq b(\varphi(t) , \sigma_t) \leq \frac{1}{\sqrt{t}}$ for $\varphi(t) \geq 2 \sqrt{t},$
by dominated convergence, if $\varphi_n$ converges towards $\varphi$, $\Lambda(\varphi_n)$ converges to $\Lambda(\varphi)$ pointwise. Now, since the functions involved
are increasing, the convergence holds uniformly on the compact $[0,1]$.\hfill $\Box$

\section{The lower bound}
\label{sec:lower}

In this section, we prove the large deviation lower bound, namely:
\begin{proposition} \label{lowerbound0}
 For any open set $O$ in $ C_\T([0,1]; \R)$, 
 \begin{equation} \label{eq-lowerbound0}
\liminf_{N \rightarrow \infty}\frac{1}{N} \ln \PP(\lambda_1^{\theta, N} \in O) \geq - \inf_{\varphi \in O}I_\theta(\varphi).
\end{equation}
\end{proposition}
 For any $\varphi \in C_\T([0,1]; \R),$ any $\delta >0,$
$B(\varphi, \delta)$ will denote the ball centered at $\varphi$ with radius
$\delta$ with respect to the uniform metric, that is the subset of  $C([0,1]; \R)$ of functions $\psi$ such that 
$ \displaystyle \sup_{t\in [0,1]} |\psi(t) - \varphi(t)| <\delta.$ 
 To prove Proposition \ref{lowerbound0}, it is enough to show
\begin{proposition} \label{lowerbound}
%Let $\T \ge 0$ be fixed.
 \begin{equation} \label{eq-lowerbound}
\lim_{\delta \downarrow 0}\liminf_{N \rightarrow \infty}\frac{1}{N} \ln \PP\left(\lambda_1^{\theta, N} \in B(\varphi, \delta)\right) \geq - I_\theta(\varphi)
\end{equation}
for any $\varphi$ belonging to a well chosen  subclass ${\mathcal H}_\theta$  of $C_\T([0,1]; \R)$ satisfying
\begin{equation} \label{densiteH}
 \inf_{\varphi \in O \cap {\mathcal H}_\theta} I_\theta(\varphi) = \inf_{\varphi \in O} I_\theta(\varphi)
 \end{equation}
for any open set $O$.
\end{proposition}
To introduce the subclass ${\mathcal H}_\theta$, we need a few more
notations.

\smallskip\noindent
 For $\varphi$ such that $I_\theta(\varphi) < \infty$, we recall from Section \ref{sec:rf} that
\begin{equation} \label{defkphi}
k_\varphi(s) := \dot{\varphi}(s) - b(\varphi(s), \sigma_s)  = \dot{\varphi}(s) - \frac{1}{2s}(\varphi(s) - \sqrt{\varphi^2(s) -4s})
\end{equation}
and that
$$ I_\T(\varphi) = \frac{1}{2}\int_0^1 k_\varphi^2(s) ds = \frac{1}{2} \|k_\varphi \|_2^2.$$
%where $\|  \cdot  \|_2$ is the usual norm in $L^2([0,1]).$
\noindent
We define
\begin{equation} \label{defH}
 \begin{array}{l}
{\mathcal H}_\theta = \{ \varphi \in C_\theta([0,1]; \R); \varphi(t)>2\sqrt{t} \  \forall t \in [0,1]; \ k_\varphi \ \rm{smooth}\} \qquad \;  \rm{for } \;  \theta>0, \\
\mbox{}\\
{\mathcal H}_0 = \left\{ \varphi \in C_0([0,1]; \R); \exists t_0 >0, \begin{array}{ll}
\varphi(t)=2\sqrt{t} & \rm{ for } \,  t \leq t_0\\
 \varphi(t)>2\sqrt{t} & \rm{ for } \, t>t_0 \  \end{array} ; k_\varphi \ \rm{ smooth } \right\}, 
\end{array} 
\end{equation}
where smooth  means infinitely differentiable on $[0,1].$
 For $\varphi \in {\mathcal H}_0,$ we denote by $t_0(\varphi):= \sup\{t; \varphi(t) = 2\sqrt t\}$ the corresponding threshold.\smallskip

\noindent
Eq.\eqref{densiteH} will be proven in Lemma \ref{lemapprox}
%In Lemma \ref{lemapprox} will be proven  that ${\mathcal H}_\theta$ is dense  in $\{ \varphi \in C_\T([0,1]); \, I_\T(\varphi) <\infty \}$ and \eqref{densiteH}
 after some preliminary considerations in the next subsection. Eq.\eqref{eq-lowerbound} is obtained in Section 4.4 for $\T >0$ and in Section 4.5 for $\T =0$.\smallskip

\subsection{Some properties of the functions with finite entropy when $\T=0$}

As will be seen further, the proof  that ${\mathcal H}_\theta$ is dense will be quite straightforward in the case when
$\T>0$ but more delicate when $\T=0.$ In this latter case, we first need to understand some features of the functions with finite entropy
that we gather here.

We need the following notations :
for any $\varphi$ such that $\varphi(s) \geq 2\sqrt{s}, \,\,  \forall s \in [ 0,1],  $ we define $x_{\varphi}$ by
\begin{equation} \label{defx}
 x_{\varphi}(s) = \frac{\varphi(s)+ \sqrt{\varphi^2(s) -4s}}{2}, \quad  \forall s \in [ 0,1]
\end{equation}
so that
$\varphi$ and $k_\varphi$ can be reexpressed in terms of $x_{\varphi}.$ More precisely, $\forall s \in ( 0,1],$
\begin{equation} \label{phifctx}
\varphi(s)  = x_{\varphi}(s) + \frac{s}{x_{\varphi}(s)} 
\end{equation}
and 
\begin{equation} 
\label{kaphi} 
k_\varphi (s) = 2\dot{x}_{\varphi}(s) \left( 1 - \frac{s}{x_{\varphi}^2(s)}\right). 
\end{equation}

\noindent
The following lemma gives the behaviour of $\varphi$ near $0.$

\begin{lemma} ($\theta = 0$) \label{lem.comp.en0}
Let $\varphi \in C_0([0,1)]$ satisfy $I_0(\varphi) < \infty$. Then,
$$
\lim_{t \rightarrow 0} \frac{\varphi(t)}{\sqrt{t}} = 2.
$$
\end{lemma}

 \noindent{\bf Proof:} Set 
$$ I^t(\varphi) = \int_0^t \left(\dot{\varphi}(s) - \frac{1}{2s}\left(\varphi(s) - \sqrt{\varphi^2(s) -4s}\right) \right)^2 ds.$$
Then, from the finiteness of $I_0(\varphi)$, $\lim_{t \rightarrow 0} I^t(\varphi) = 0$.
 From Cauchy-Schwarz inequality,
$$ \left| \int_0^t \dot{\varphi}(s) - \frac{1}{2s}\left(\varphi(s) - \sqrt{\varphi^2(s) -4s}\ \right) ds \right| \leq \sqrt{t} (I^t(\varphi))^{1/2}$$
and
$$ \left| \frac{\varphi(t)}{\sqrt{t}} - \frac{1}{\sqrt{t}} \int_0^t \frac{1}{2s}\left(\varphi(s) - \sqrt{\varphi^2(s) -4s}\right) ds \right| \leq (I^t(\varphi))^{1/2}.$$
Now, we have, using \eqref{phifctx},
$$ 0 \leq  \int_0^t \frac{1}{2s}\left(\varphi(s) - \sqrt{\varphi^2(s) -4s}\right) ds = \int_0^t \frac{ds}{x_\varphi(s)} \leq \int_0^t \frac{ds}{\sqrt{s}}  = 2 \sqrt{t}.$$
Thus, on one hand, $\displaystyle \frac{\varphi(t)}{\sqrt{t}} \geq 2$, whereas $ \displaystyle 0 \leq \frac{1}{\sqrt{t}} \int_0^t \frac{1}{2s}\left(\varphi(s) - \sqrt{\varphi^2(s) -4s}\right) ds \leq 2$ and the difference of the two terms tends to 0 as $t$ tends to 0. It follows that:
$$ \lim_{t \rightarrow 0} \frac{\varphi(t)}{\sqrt{t}} =  \lim_{t \rightarrow 0}\frac{1}{\sqrt{t}}\int_0^t \frac{1}{2s}\left(\varphi(s) - \sqrt{\varphi^2(s) -4s}\right) ds = 2. $$
\hfill $\Box$

\noindent
The following lemma will be useful in the proof of the lower bound itself.
% When $\theta = 0$, we shall need  the following
 \begin{lemma} \label{lempos} ($\theta = 0$)
  Let $\varphi \in {\mathcal H}_0$. 
  Then $k_\varphi$ is positive in a right neighborhood of $t_0(\varphi)$. 
  \end{lemma}

\noindent
{\bf Proof of Lemma \ref{lempos}:} 
 For $\varphi \in {\mathcal H}_0$, $k_\varphi \equiv 0$ on $[0, t_0(\varphi)]$.

\noindent
 Since $\varphi(s) > 2\sqrt{s}$ for $s> t_0(\varphi)$, we have that $x_\varphi(s) > \sqrt{s}$,  for $s > t_0(\varphi).$
As  $\dot{x}_\varphi(t_0(\varphi))= \frac{1}{2\sqrt{t_0(\varphi)}} >0$ and $\dot x_\varphi$ is continuous (as $\varphi$ is smooth),
$\dot{x}_\varphi(s) >0 $ in a neighborhood of $t_0(\varphi)$ and thus, from \eqref{kaphi}, $k_\varphi(s) > 0$ for $t_0(\varphi)<s < t_0(\varphi) + \varepsilon$ 
for some $\varepsilon >0$. \hfill $\Box$

\subsection{Denseness of ${\mathcal H}_\theta$}

The goal of this subsection is to establish the following lemma

\begin{lemma} \label{lemapprox}
Let $\varphi \in C_\T([0,1])$ satisfying $I_\T(\varphi) < \infty.$ 
 There exists a sequence $(\varphi_p)_{p \in \mathbb N^*}$ of functions in $ {\mathcal H}_\theta$ such that,
as $p$ goes to infinity,

%\item[-] $\dot{\varphi}_p(s) - b(\varphi_p(s), \sigma_s)$ is absolutely continuous
 - $\varphi_p$ converges to $\varphi$ in $C_\T([0,1], \R)$
 
 - $I_\T(\varphi_p)$ converges to $I_\T(\varphi)$.
Therefore, \eqref{densiteH} holds.
\end{lemma}

\subsubsection{Proof of Lemma \ref{lemapprox} when $\T>0$}

Let $\varphi$ such that $I_\T(\varphi) < \infty$. 
As $\varphi(0) = \T >0$ and $\varphi$ is continuous, there exists $t_1>0$
such that for any $t \in [0, t_1],$ $\varphi(t) > 2\sqrt{t}.$

 For any $p \in \mathbb N^*,$ we define
\[ \chi_p(t) =\left\{ 
\begin{array}{ll}
 \varphi(t) & \textrm{if } t\leq t_1,\\
 \varphi(t) + (t-t_1) & \textrm{if }  t_1\leq t \leq t_1 + \frac 1 p,\\
 \varphi(t) + \frac 1 p & \textrm{if }   t \geq t_1 + \frac 1 p.\\
\end{array}
\right.
\]

It is easy to check that $\chi_p$ is continuous and  for $p$ large enough, for any $t \in [0, 1],$ $\chi_p(t) > 2\sqrt{t}$ 
and $I_\T(\chi_p) < \infty.$ Moreover, as $p$ goes to infinity, $\chi_p$ converges to 
$\varphi$ in the uniform norm and $k_{\chi_p}$ converges to $k_\varphi$ in $\mathbb L^2([0,1]).$

It is now enough to check that  $\chi_p$ can be approximated by a sequence of functions in ${\mathcal H}_\T.$
As we know that for any $t \in [0, 1],$ $\chi_p(t) > 2\sqrt{t},$ we have that 
$\inf_{s \in (0,1]} \left( 1 - \frac{s}{x_{\chi_p}^2(s)}\right) >0.$
As $k_{\chi_p} \in \mathbb L^2([0,1]),$ from \eqref{kaphi} we get that so does $\dot{x}_{\chi_p}.$ It can be approximated by a sequence of smooth functions $\dot{x}_{p,q}$. Set $x_{p,q}(t) = \theta + \int_0^t \dot{x}_{p,q}(s) ds.$  The corresponding $\chi_{p,q}$ is defined by 
$$ \chi_{p,q}(s)  = x_{p,q}(s) + \frac{s}{x_{p,q}(s)}$$ and 
$$k_{p,q}(s) =2 \dot{x}_{p,q}(s) \left( 1 - \frac{s}{x^2_{p,q}(s)}\right),$$
so that $k_{p,q}$ is smooth. For $q$  large enough, for any $s \in [0,1],$ $ x_{p,q}(s) > \sqrt{s}$, so that $\chi_{p,q}(s) > 2 \sqrt{s}.$ 

Moreover, as $q$ grows to infinity, the sequence $x_{p,q}$ converges towards $x_{\chi_p}$ in uniform norm on $[0,1]$ so that 
 $\chi_{p,q}$ converges towards $\chi_p$ in the same sense and $k_{p,q}$
converges to $k_{\chi_p}$ in $\mathbb L^2([0,1]).$\smallskip

\noindent
To conclude the proof of the lemma, it is enough to notice that one can find an increasing function $\psi$
from $\mathbb N$ to $\mathbb N$ such that $\varphi_p :=\chi_{p,\psi(p)} \in {\mathcal H}_\T$ converges towards $\varphi$ and $k_{\varphi_p} = k_{p,\psi(p)}$
converges to $k_{\varphi}$ in $\mathbb L^2([0,1]).$

\subsubsection{Proof of Lemma \ref{lemapprox} when $\T=0$}

As in the latter paragraph, we establish the proof in two steps: first, we approximate $\varphi $ by 
a sequence of functions that are equal to 
$2 \sqrt{t}$ in a neighborhood of 0 and strictly greater than $2 \sqrt{t}$ away from 0.
Next, we approximate those functions by  smooth ones.

 Let $r>0$ and define $\chi_r$ by:
 $$ \chi_r(t) = \left\{ \begin{array}{ll}
 2\sqrt{t} & t \leq y_r^2, \\
 
  \displaystyle  y_r + \frac{t}{y_r} & y_r^2 \leq t \leq r,\\
   \varphi(t) + r (t-r), & t \geq r
  \end{array} \right. $$
  with $ y_r =  \frac{\varphi(r) -  \sqrt{\varphi^2(r) -4r}}{2} \leq \sqrt{r}$ so that $\chi_r$ is continuous.
 $$ \Vert \varphi - \chi_r \Vert \leq \sup_{s \leq r} \vert \varphi(s) - \chi_r(s) \vert  \vee r \leq 2 \sup_{s \leq r} \vert \varphi(s) - 2 \sqrt{s} \vert \vee r \rightarrow_{r \rightarrow 0} 0$$
 using Lemma \ref{lem.comp.en0}.
 It remains to show that $I_0(\varphi) - I_0(\chi_r)$ tends to 0.
If we set 
$J_r(f) = \int_{0}^r k_f^2(s)ds$ and $J^\prime_r(f) = \int_{r}^1 k_f^2(s)ds,$ we get 
 $$I_0(\varphi) - I_0(\chi_r) = (J_r(\varphi) - J_r(\chi_r))+(J^\prime_r(\varphi) - J^\prime_r(\chi_r)) $$
 with $J_r(\varphi) \rightarrow 0$ as $r \rightarrow 0$.
 $$J_r(\chi_r) = \int_{y_r^2}^r \left(\frac{1}{y_r} -\frac{y_r}{s}\right)^2 ds = \frac{r}{y_r^2} - \frac{y_r^2}{r}  -  2 \ln\left( \frac{r}{y_r^2}\right) \rightarrow_{r \rightarrow 0} 0$$
 since $y_r/\sqrt{r}$ tends to 1 thanks to Lemma \ref{lem.comp.en0}.

On the other hand, if we define $h_{\varphi,r}$ by 
$$h_{\varphi,r}(t) := r- \frac{1}{2t} \left( r(t-r) + \sqrt{\varphi^2(t)-4t} - \sqrt{(\varphi(t)+r(t-r))^2-4t} \right), $$
then, by Cauchy-Schwarz inequality, we have  
$$ \vert J^\prime_r(\varphi) - J^\prime_r(\chi_r) \vert^2 \leq  
\int_r^1 (2 k_\varphi(s) + h_{\varphi,r}(t) )^2 dt\int_r^1 ( h_{\varphi,r}(t) )^2 dt. $$
Therefore it is enough to show that
$\int_r^1 (h_{\varphi,r}(t) )^2 dt $  goes to zero as $r$ goes to zero.

To show that, we notice that, for $t \in [r,1],$
\begin{equation}\label{rtr}
  \left\vert \frac{r(t-r)}{2t}  \right\vert\leq \frac{r}{2}.
  %+ \frac{r}{2} \leq r.
\end{equation}
Moreover, 
\begin{eqnarray}
\left\vert \frac{1}{2t} \left( \sqrt{\varphi^2(t)-4t} \right. \right.&- & \left.\left. \sqrt{(\varphi(t)+r(t-r))^2-4t}\right)\right\vert \\
& = &  \frac{1}{2t} \frac{2r \varphi(t)(t-r)+r^2(t-r)^2}{ \sqrt{\varphi^2(t)-4t} + \sqrt{(\varphi(t)+r(t-r))^2-4t}} \nonumber\\
& \leq &  \frac{\sqrt{2r \varphi(t)(t-r)}}{2t}  +\frac{|r(t-r)|}{2t} \nonumber\\
& \leq &  C r^{1/4}+ \frac{r}{2}, \label{J}
\end{eqnarray}
where we used that, from Lemma \ref{lem.comp.en0}, $t \mapsto \frac{\varphi(t)}{\sqrt t}$ is bounded
on $[0,1]$ by a constant $C.$
Putting \eqref{rtr} and \eqref{J} together, we get 
$\int_r^1 (h_{\varphi,r}(t) )^2 dt $  goes to zero as $r$ goes to zero.

 Now $\dot{\chi}_r(s) = \frac{1}{\sqrt{s}}$  on $[0, y_r^2]$ and $\dot{\chi}_r \in \mathbb L^2([y_r^2, 1])$ since $k_{\chi_r} \in \mathbb L^2([y_r^2, 1]).$  
 For any $r>0,$ there exists a sequence of function $\dot{\chi}_{r,q}$ smooth on $]0,1]$ such that $\dot{\chi}_{r,q} (s) = \frac{1}{\sqrt{s}}$
  on $[0, y_r^2]$ and $\dot{\chi}_{r,q}$ tends to $\dot{\chi}_r$ in $\mathbb L^2([y_r^2/2, 1])$. 
 Setting $\chi_{r,q}(t) = \int_0^t \dot{\chi}_{r,q}(s) ds$, we have: 
\begin{itemize}
 \item[-]  $\chi_{r,q}$ tends to $\chi_r$ in uniform norm.
\item[-] $k_{\chi_{r,q}}$ is smooth.
\item[-]  $k_{\chi_{r,q}}$ converges to $k_{\chi_r}$ in $\mathbb L^2([0,1])$ so that $I_0(\chi_{r,q})$ converges to $I_0(\chi_r)$.
\end{itemize}

Putting everything together, we conclude that there exists an increasing function $\psi$ such that
 the sequence of functions $\varphi_p = \chi_{p,\psi(p)}$ satisfies the requirements of  Lemma \ref{lemapprox}. \hfill $\Box$

%\vspace{.3cm}

%%%%%%%%%%%%%%%%%%%%%%%%%%%%%%%%%%%%%%%

\subsection{Almost sure convergence of $\lambda_1^{\theta,N}$ and $\nu_N$ under $\PP^{k_\varphi}$}

The strategy of the proof of the lower bound will be classical : we will make a proper change of measure
so that the function $\varphi$ becomes ``typical'' under the new measure. We will therefore need to study 
more precisely the behavior of $\lambda_1^{\theta,N}$ under the new measure $\PP^{k_\varphi}.$ 

To be more precise, 
for $h \in \mathbb L^2([0,1]),$ we define 
 the exponential martingale $M^h$ such that for any $t \in [0,1],$
\begin{equation}\label{defimart}
  M_t^h = \exp \left[ N \left( \int_0^t h(s) \frac{1}{\sqrt{N}}
  dB_1(s) - \frac{1}{2} \int_0^t h^2(s) ds\right) \right],
\end{equation}
where $B_1$ is the standard Brownian motion appearing in the SDE
for $\lambda_1^{\theta,N}$ (see \eqref{EDSvp}). We denote by $(\mathcal F_t)_{t\ge 0}$ its canonical filtration.

We now introduce 
 $\PP^{ k_\varphi}$ the probability defined by $\PP^{ k_\varphi}
:= M^{ k_\varphi}_1 \sharp \PP,$ meaning that  for any $t \le 1,$ the Radon-Nikodym derivative of $\PP^{ k_\varphi}$
with respect to $\PP$ on $\mathcal F_t$ is given by  $M^{ k_\varphi}_t$ and we also denote by $\E^{ k_\varphi}$ the expectation
under $\PP^{ k_\varphi}.$
Recall that  $\nu_N$ is the empirical distribution of all but the largest eigenvalues
defined in the introduction.
\noindent
 For any $r>0, \alpha >0,$ we also define 
$$\mathbb B_r(\sigma, \alpha) := \mathbb B(\sigma, \alpha) \bigcap \{ \mu \in C([0,1],
{\mathcal P}(\R)); \, \forall s, \,  \textrm{supp}(\mu_s) \subset ]-\infty, 2\sqrt{s}+r] \}.$$ 

The goal of this subsection will be to show
\begin{proposition} \label{propPh}
 For any $r>0, \delta >0, \alpha>0$ and $\varphi \in {\mathcal H}_\T,$

\smallskip\noindent
 $$\PP^{k_\varphi}(\lambda_1^{\T,N} \in B(\varphi, \delta);\nu_N \in \mathbb B_r(\sigma, \alpha))) 
\xrightarrow[N \rightarrow \infty]{}1.$$% tends to 1 as $N$ goes to $\infty$.
\end{proposition}
The proof of the proposition relies on some lemmata.
\begin{lemma} \label{lem1}
Under $\PP^{k_\varphi}$, $\mu_N$ and $\nu_N$ converge as $N$ goes to infinity to the semicircular process $\sigma$.
\end{lemma}

\noindent
{\bf Proof:}  It is well known that $\mu_N$ is exponentially tight in scale $N^2$, under $\PP$ (see \cite{CDG}, \cite[Chap. 12]{G}). Let $A \in C([0,1], {\mathcal P}(\R))$, then
\begin{eqnarray}
\PP^{k_\varphi}( \mu_N \in A) &=& \E\left( M_1^{k_\varphi} {\bf 1}_{\mu_N \in A}\right) \nonumber \\
&\leq & \E( (M_1^{k_\varphi})^2)^{1/2} (\PP( \mu_N \in A))^{1/2} \nonumber \\
&=& \exp\left(\frac{N}{2} \int_0^1 k_\varphi^2(s) ds\right) \PP( \mu_N \in A)^{1/2}. \label{tension}
\end{eqnarray}
 From the exponential tighness of $\mu_N$ under $\PP$, there exists a compact $K_L$ in $C([0,1], {\mathcal P}(\R))$ such that 
$$ \PP(\mu_N \in K_L^c) \leq \exp( - N^2 L).$$
Therefore, from \eqref{tension}
$$ \PP^{k_\varphi} (\mu_N \in K_L^c) \leq \exp( - N^2 L/4)$$
for N large enough. This proves the exponential tightness of $\mu_N$ under $\PP^{k_\varphi}$ and thus its a.s. pre-compactness  in $ C([0,1], {\mathcal P}(\R))$. It remains to prove the uniqueness of any limit point. 

 From Girsanov's theorem, we have that under $\PP^{k_\varphi}\!$, the process $(\!\lambda_i^{\T,N}\!(t)\!)_{t \leq 1, i= 1, \ldots N}$
satisfies the system of stochastic differential equations:
\begin{equation} \label{EDSvp-sous-h}
\left \{ \begin{array}{l} \displaystyle d\lambda_i(t) = \frac{1}{\sqrt{N}}d\beta_i(t) + \frac{1}{N} \sum_{ j \not= i} \frac{1}{\lambda_i(t) - \lambda_j(t)} dt , \qquad   \qquad \; i = 2, \ldots, N \\
\displaystyle d\lambda_1(t) = \frac{1}{\sqrt{N}}d\beta_1(t) + k_\varphi(t) dt + \frac{1}{N} \sum_{ j \not= 1} \frac{1}{\lambda_1(t) - \lambda_j(t)} dt
\end{array} \right.
\end{equation}
where $(\beta_i)_{1\le i\le N}$ are independent Brownian motions under $\PP^{k_\varphi}$.
The proof of  the uniqueness of any limit point  follows the same proof as in \cite[Theorem 1]{RS} (see also \cite[Chap. 12]{G}, \cite{C}):  let $f \in C^2_b(\R)$, using It\^o's formula and \eqref{EDSvp-sous-h}, we obtain a stochastic differential equation for $\langle \mu_N(t),f \rangle$ with a diffusion coefficient tending to 0 as $n$ tends to $\infty$. Then,  any limit point $\mu_t$ satisfies a deterministic evolution equation (the term in $k_\varphi$ disappears in the limit)
$$ \langle \mu_t , f \rangle =  \int f(x) d\mu_t(x) = \int f(x) d\mu_0(x) +\frac{1}{2}  \int_0^t \int \frac{f'(x) - f'(y)}{x-y} d\mu_s(x) d\mu_s(y) ds$$
 for which uniqueness holds. When $ \mu_0= \delta_0$  as in our setting, $\mu_t$ is the semicircular law $\sigma_t$. Therefore $\mu_N$ converges a.s. to 
the semicircle process $\sigma$.  Since $d(\mu_N, \nu_N) \leq \frac{2}{N}$, the same convergence holds for $\nu_N$.\hfill $\Box$

 \begin{lemma} For any $r>0, \alpha>0,$

\smallskip\noindent
$$\PP^{k_\varphi}( \nu_N \in \mathbb B_r( \sigma, \alpha)) \xrightarrow[N \rightarrow \infty]{} 1.$$ %converges to $1$ as $N$ goes to $\infty$.
 \end{lemma}
 
\noindent{\bf Proof:} Since we already know the convergence of $\nu_N$ towards $\sigma$ under $\PP^{k_\varphi}$, it is enough to prove that 
 under $\PP^{k_\varphi}$, 
\begin{equation} \label{lambda2}
\limsup_{N\rightarrow \infty}\lambda_2(t)
 \leq 2 \sqrt{t}.
 \end{equation}
We define  $(\lambda_i^{(\varepsilon)}(t), i= 1, \ldots N)$
the strong solution of the system of SDE  \eqref{EDSvp-sous-h} with
 initial conditions $\lambda_1^{(\varepsilon)}(0) = \T + \varepsilon$
 and
$\lambda_i^{(\varepsilon)}(0) = \frac \varepsilon i,$ for $i =
 2,\ldots, N,$ so that in particular
 $\lambda_2(t)=\lambda_2^{(0)}(t).$

We also introduce 
$(\overline{\lambda_i^{(\varepsilon)}}, i= 2, \ldots, N)$ the strong
 solution of the system of SDE:
 \begin{equation} \label{EDS(N-1)}
 d\overline{\lambda_i^{(\varepsilon)}}(t) = \frac{1}{\sqrt{N}}d\beta_i(t) +
 \frac{1}{N} \sum_{j=2,  j \not= i}^N
 \frac{1}{\overline{\lambda_i^{(\varepsilon)}}(t) -
 \overline{\lambda_j^{(\varepsilon)}}(t)} dt , \qquad   \qquad \; i = 2,
 \ldots, N,
 \end{equation}
with
 initial conditions $\overline{\lambda_i^{(\varepsilon)}}(0) = \frac \varepsilon i,$ for $i =
 2,\ldots, N.$

The process $(\overline{\lambda_i^{(0)}}, i= 2, \ldots, N)$ is distributed as the eigenvalues of $\sqrt{\!\frac{N-1}{N}}\!H_{N-1}\!(t)$ where $H_{N-1}$ is a standard Hermitian Brownian motion of size $N-1$. Therefore, $\lim_{N \rightarrow \infty}
\overline{\lambda_2^{(0)}}(t) = 2 \sqrt{t}$ a.s. 

Our goal is now to compare  $\lambda_2^{(0)}(t)$ with $\overline{\lambda_2^{(0)}}(t).$ 

The first step is to show that for any $\varepsilon >0$ fixed, $N$
 fixed, for all $t \in [0,1],$  
$\lambda_2^{(\varepsilon)}(t) \leq
 \overline{\lambda^{(\varepsilon)}_2}(t).$
In fact, we will show that for any $2\le i \le N,$ we have 
$\lambda_i^{(\varepsilon)}(t) \leq
 \overline{\lambda^{(\varepsilon)}_i}(t),$ for all $t \in [0,1].$ 

Let $R > 0$ large enough so that $\frac 1 R <
 \frac{\varepsilon}{N^2}$
and 
%%%%%%%%%%
%% Overfull 9pt
%%%%%%%%%%
% $$T_R = \inf \left\{t \geq 0, \forall i,j = 2, \ldots, N, i \neq j,\,\,
% \left|  \lambda^{(\varepsilon)}_i(t) -
%  \lambda^{(\varepsilon)}_j(t)\right| \vee
%  \left|\overline{\lambda^{(\varepsilon)}_i}(t)
%  -\overline{\lambda^{(\varepsilon)}_j}(t)\right| 
% \leq \frac 1 R\right\}.$$
\begin{eqnarray*}
 T_R &=&         \inf \left\{t \geq 0, \forall i,j = 2, \ldots, N, i \neq j, \phantom{\frac 1 R}\right.\\
      &&\phantom{\inf \left\{\right.}\left.\left|  \lambda^{(\varepsilon)}_i(t) -
 \lambda^{(\varepsilon)}_j(t)\right| \vee
 \left|\overline{\lambda^{(\varepsilon)}_i}(t)
 -\overline{\lambda^{(\varepsilon)}_j}(t)\right| 
\leq \frac 1 R\right\}.
\end{eqnarray*}
 For any $t \ge 0$ and   $2\le i \le N,$ we can write
$$d\lambda_i^{(\varepsilon)}(t) = \frac{1}{\sqrt{N}}d\beta_i(t) + f_i(\boldsymbol{\lambda}^{(\varepsilon)}(t))dt,$$
and
$$d\overline{\lambda_i^{(\varepsilon)}}(t) = \frac{1}{\sqrt{N}}d\beta_i(t) + g_i(\overline{\boldsymbol\lambda^{(\varepsilon)}}(t))dt,$$
with $\boldsymbol{\lambda}^{(\varepsilon)}(t):= (\lambda_1^{(\varepsilon)}(t), \ldots,\lambda_N^{(\varepsilon)}(t) ),$
$\overline{\boldsymbol\lambda^{(\varepsilon)}}(t):= (\overline{\lambda_2^{(\varepsilon)}}(t), \ldots,\overline{\lambda_N^{(\varepsilon)}}(t) )$
and $f_i:\R^N \rightarrow \R$ and $g_i:\R^{N-1} \rightarrow \R$ are such that
$$ f_i(x_1,\ldots,x_N)= \frac{1}{N} \frac{1}{x_i-x_1} + g_i(x_2,\ldots,x_N).$$
We also denote by $\widetilde{\boldsymbol\lambda^{(\varepsilon)}}(t):= (\lambda_2^{(\varepsilon)}(t), \ldots,\lambda_N^{(\varepsilon)}(t) ).$
Then
$$ d\left(\lambda_i^{(\varepsilon)}-\overline{\lambda_i^{(\varepsilon)}}\right)(t)= (f_i(\boldsymbol{\lambda}^{(\varepsilon)}(t))
-g_i(\overline{\boldsymbol\lambda^{(\varepsilon)}}(t))dt \le 
(g_i(\widetilde{\boldsymbol{\lambda}^{(\varepsilon)}}(t))
-g_i(\overline{\boldsymbol\lambda^{(\varepsilon)}}(t)))dt,$$
as $\lambda_i^{(\varepsilon)}(t)- \lambda_1^{(\varepsilon)}(t)\le 0.$\smallskip

We denote by $x^+= \max(0,x)$ so that 
$$ d\left(\lambda_i^{(\varepsilon)}-\overline{\lambda_i^{(\varepsilon)}}\right)^+(t)
\le \mathbf{1}_{\overline{\lambda_i^{(\varepsilon)}}(t)\le \lambda_i^{(\varepsilon)}(t)}
(g_i(\widetilde{\boldsymbol{\lambda}^{(\varepsilon)}}(t))
-g_i(\overline{\boldsymbol\lambda^{(\varepsilon)}}(t)))dt.$$

Now
$$g_i(\widetilde{\boldsymbol{\lambda}^{(\varepsilon)}}(t))
-g_i(\overline{\boldsymbol\lambda^{(\varepsilon)}}(t)) 
= \frac{1}{N} \sum_{\substack{
                     k \neq i\\
                     k \ge 2}
} \frac{(\overline{\lambda_i^{(\varepsilon)}}-\lambda_i^{(\varepsilon)})-(\overline{\lambda_k^{(\varepsilon)}}-\lambda_k^{(\varepsilon)})}
{(\lambda_i^{(\varepsilon)}-\lambda_k^{(\varepsilon)})(\overline{\lambda_i^{(\varepsilon)}}-\overline{\lambda_k^{(\varepsilon)}})}(t)$$
As the eigenvalues are ordered, the denominator is always positive so that, for all $i \ge 2,$
$$ \mathbf{1}_{\overline{\lambda_i^{(\varepsilon)}}(t)\le \lambda_i^{(\varepsilon)}(t)}
\frac{(\overline{\lambda_i^{(\varepsilon)}}-\lambda_i^{(\varepsilon)})}
{(\lambda_i^{(\varepsilon)}-\lambda_k^{(\varepsilon)})(\overline{\lambda_i^{(\varepsilon)}}-\overline{\lambda_k^{(\varepsilon)}})}(t) \le 0.
$$
On the other hand, for all $k \neq i,$ for $t \in [0, T_R],$
$$\frac{-(\overline{\lambda_k^{(\varepsilon)}}-\lambda_k^{(\varepsilon)})}
{(\lambda_i^{(\varepsilon)}-\lambda_k^{(\varepsilon)})(\overline{\lambda_i^{(\varepsilon)}}-\overline{\lambda_k^{(\varepsilon)}})}(t) \le
R^2 (\lambda_k^{(\varepsilon)}-\overline{\lambda_k^{(\varepsilon)}})^+(t),$$
so that on $[0, T_R],$
$$ d\left(\lambda_i^{(\varepsilon)}-\overline{\lambda_i^{(\varepsilon)}}\right)^+(t) \le \frac{R^2}{N} \sum_{k \neq i}
(\lambda_k^{(\varepsilon)}-\overline{\lambda_k^{(\varepsilon)}})^+(t)dt,$$
and if we sum over the index $i,$
$$ d\left(\sum_{i=2}^N(\lambda_i^{(\varepsilon)}-\overline{\lambda_i^{(\varepsilon)}})\right)^+(t) \le \frac{R^2 (N-1)}{N} \sum_{k =2}^N
(\lambda_k^{(\varepsilon)}-\overline{\lambda_k^{(\varepsilon)}})^+(t)dt.$$
By Gronwall lemma, we get that on $[0, T_R],$
$$\sum_{i=2}^N(\lambda_i^{(\varepsilon)}-\overline{\lambda_i^{(\varepsilon)}})(t)=0,$$
which means that for all $i \ge 2,$ 
$\lambda_i^{(\varepsilon)}(t) \leq
 \overline{\lambda^{(\varepsilon)}_i}(t).$ 

\smallskip
Moreover, from \cite{CL}, we know that $T_R$ goes to infinity 
as $R$ goes to infinity. In particular, if we choose $R$ large enough
 for $T_R$ to be larger than 1, our inequalities hold for any $t \in
 [0,1].$

\smallskip\noindent
Now, the solutions of \eqref{EDS(N-1)}, resp. \eqref{EDSvp-sous-h}, are continuous with respect to the initial condition (see \cite{Ce}); thus, letting $\varepsilon \rightarrow 0$, we obtain 
$\lambda_2^{(0)}(t) \leq
 \overline{\lambda^{(0)}_2}(t)$ a.s.
Putting everything together, we have that 
$\limsup_{N\rightarrow \infty}\lambda_2(t)
 \leq 2 \sqrt{t}.$ \hfill
$\Box$

 \begin{lemma}
Let $\varphi \in \mathcal H_\T$ . 
\begin{itemize}
\item[1)] For $\theta >0,$ the differential equation  $ dy(t) = (k_\varphi(t) + b(y(t), \sigma_t)) dt$ on $[0,1]$ with initial value $y(0)=\theta$ admits a unique solution larger than 
 $2 \sqrt{t}$, namely $\varphi$. 
\item[2)]   For $\theta =0,$ the differential equation  $ dy(t) = (k_\varphi(t) + b(y(t), \sigma_t)) dt$  on $[0,1]$ with  value 
 $y(t_0(\varphi))=2 \sqrt{t_0(\varphi)}$ at time $t_0(\varphi)$  admits a unique solution larger than 
 $2 \sqrt{t}$, namely $\varphi$. 
\end{itemize}
 \end{lemma}

  \noindent{\bf Proof:} 
 Let
 us first check that in both cases there is a unique solution larger than $2\sqrt t$. We recall that for $x \geq  2\sqrt t,$
 $$b(x, \sigma_t )= \frac{1}{2\pi t}\  \int \frac{1}{x-y} \sqrt{ 4t - y^2} dy = \frac{1}{2t} ( x - \sqrt{ x^2 - 4t}) .$$ It is easy to see that $x \mapsto b(x, \sigma_t)$ is decreasing on $ [2\sqrt{t}, \infty[$.
 Let $x, y$ two solutions of $ dy(t) = (k_\varphi(t) + b(y(t),
 \sigma_t)) dt$ such that for any $t \in [0,1],$ $ x(t), y(t) \geq 2 \sqrt{t}$. Then,
 $$ (x(t) - y(t))^2 = 2 \int_0^t (x(s) - y(s)) (b(x(s), \sigma_s) - b(y(s), \sigma_s)) ds \leq 0.$$
In the first case, it is very easy to check that $\varphi$ is a solution.
In the second case, notice that for $t \leq t_0(\varphi),$ $k_\varphi(t) = 0$ and we know that
$t \mapsto 2\sqrt t$ is a solution of  $ dy(t) =  b(y(t), \sigma_t) dt$ with initial condition $y(0)=0.$ $\hfill \Box$

The last lemma to complete the proof of Proposition \ref{propPh} is the following
\begin{lemma} \label{lemder}
 For any $\T \geq 0$ and $\varphi \in \mathcal H_\T,$ under $\PP^{k_\varphi}$, the process $\lambda_1^{\T,N}$ converges a.s. to $\varphi$.
\end{lemma}

 \noindent
\textbf{Proof :} 
 As in Lemma \ref{lem1}, from the exponential tightness in scale $N$ of $\lambda_1$ under $\PP$, we deduce the exponential tighness of $\lambda_1$ under $\PP^{k_\varphi}$ and the a.s. pre-compactness  of $\lambda_1$.
 Let $x(t)$ be a limit point. There exists $f : \mathbb N \rightarrow
 \mathbb N$ strictly increasing such that $\lambda_1^{\T,f(N)}(t)$ converge
 to $x(t).$ In the sequel we omit the superscript $\T,f(N).$

\smallskip\noindent
The crucial step of the proof, which is similar for any value of $\T$ is to show that $x(t) \geq \varphi(t).$

\smallskip
 From the a.s. convergence of $(\mu_N)_t$ towards $\sigma_t$ and using that $\sigma_t( [2 \sqrt{t} - \varepsilon, 2\sqrt{t}]) >0$, it follows that $\liminf_N \lambda_1(t) \geq 2 \sqrt{t}$ and thus
 $x(t) \geq 2\sqrt{t}$.
 From It\^o's formula, we get
\begin{eqnarray*} 
((\varphi(t) - \lambda_1(t))^+)^2 & = & - \frac{ 2}{ \sqrt N} (\varphi(t)
 - \lambda_1(t))^+ d\beta_1(t)\\
& & + 2 \int_0^1 (\varphi(t) - \lambda_1(t))^+ [b(\varphi(t), \sigma_t) -
 b_N(\lambda_1(t), (\nu_N)_t)] dt\\ 
& & + \frac 1 N{\bf 1}_{\varphi(t) -
 \lambda_1(t)\geq 0} dt
\end{eqnarray*}
where $b_N = \frac{N-1}{N}b$.

The first and last term converge to zero and 
we decompose the second term in three  $2 (A^1(t) + A^2(t) + A^3(t) )$ where
\begin{multline*}
 A^1(t) = \int_0^1 (\varphi(t) - \lambda_1(t))^+ [b(\varphi(t), \sigma_t) -
 b_N(\varphi(t), \sigma_t)] dt \\
= \frac{1}{N}  \int_0^1 (\varphi(t) - \lambda_1(t))^+ b(\varphi(t), \sigma_t) dt,
\end{multline*}
 
 $$ A^2(t) =  \frac{N-1}{N} \int_0^1 (\varphi(t) - \lambda_1(t))^+ [b(\varphi(t), \sigma_t) -
 b(\lambda_1(t), \sigma_t)] dt $$
and $$ A^3(t) =\frac{N-1}{N}  \int_0^1 (\varphi(t) - \lambda_1(t))^+ [b(\lambda_1(t), \sigma_t) -
 b(\lambda_1(t), (\nu_{N})_t)] dt. $$
 Passing to the limit, we obtain:
$$ ((\varphi(t) - x(t))^+)^2 =2  \lim_{N \rightarrow \infty} (A^1(t) +
 A^2(t) +A^3(t) )\; \rm{a.s.}.$$
We now use the continuity  on $\R$ and  the monotony on $[2 \sqrt{t}, \infty[$ of the function $x \mapsto b(x, \sigma_t)$,  the lower semicontinuity of $(x, \mu) \mapsto b(x, \mu)$ to conclude that :
$$  \lim_{N \rightarrow \infty} A^2(t) =  \int_0^1 (\varphi(t) - x(t))^+ [b(\varphi(t), \sigma_t) -
 b(x(t), \sigma_t)] dt \leq 0 \; \rm{a.s.}$$
 and
 $$ \lim_{N \rightarrow \infty} A^3(t) \leq 0 \; \rm{a.s.}.$$
We also easily get that 
 $$ \lim_{N \rightarrow \infty} A^1(t) = 0  \; \rm{a.s.}$$ since $ b(\varphi(t), \sigma_t) \in L^1([0,1])$ for $\varphi \in {\mathcal H}_\T$.
 Therefore, we obtain that  $x(t) \geq \varphi(t).$

\smallskip
In the case when $\T >0,$
 we therefore get that  $x(t)$ is
well separated from the\break support of $\sigma_t$, we can argue as before, using \eqref{lambda2},
 and obtain that\break
 $\lim_{N \rightarrow \infty} b_N(\lambda_1(t), (\nu_N)_t) = b(x(t),
 \sigma_t)$. Therefore, letting $N \rightarrow \infty$ in the equation
 of $\lambda_1$, we obtain that $x$ is a solution of the differential
 equation $ dy(t) = (k_\varphi(t) + b(y(t), \sigma_t)) dt$ and
 therefore equal to $\varphi.$

\smallskip\noindent
In the case when $\T =0,$ we have to treat first the case $t \leq t_0.$ On this interval, $k_\varphi =0,$
so that $\lambda_1(t)$ converges to $2\sqrt t,$ that is $\varphi(t).$

 For any $t >t_0,$ $x(t)$ is
well separated from the support of $\sigma_t$, and we get as before
$\displaystyle\lim_{N \rightarrow \infty} b_N(\lambda_1(t), (\nu_N)_t) = b(x(t), \sigma_t)$ so that 
 $x$ is a solution of the differential equation $ dy(t) = (k_\varphi(t) + b(y(t), \sigma_t)) dt$ with initial condition $x(t_0)=2\sqrt{t_0} $ and therefore equal to $\varphi.$ \hfill
 $\Box$

% \vspace{.3cm}
\noindent
Proposition \ref{propPh} is straightforward  from Lemmata \ref{lem1} to \ref{lemder}.

%%%%%%%%%%%%%%%%%%%%%%%%%%%%%%%%%%%%%%%
\subsection{ Lower bound for  a  non null initial condition: $\theta >0$}
%\label{lb1}
%Let $\T >0,$ $(H_N^\T(t))_{t \in [0,1]}$ the corresponding Hermitian Brownian motion  and
%$\lambda_1^\T$ the process of its largest eigenvalue.
%Let $\varphi \in C_\T([0,1], \R)$ such that  $I_\T(\varphi) < \infty.$ 
%Then 
%$$ \lim_{\delta \downarrow 0} \liminf_{N \rightarrow \infty}
%\frac{1}{N} \ln \PP(\lambda_1^\T \in B(\varphi, \delta)) \geq - I_\T(\varphi).$$
%end{prop}
We want to show Proposition \ref{lowerbound} - Eq.\eqref{eq-lowerbound}  for $\varphi \in {\mathcal H}_\theta$ under the assumption that $\theta >0.$

%Thanks to Lemma \ref{lemapprox} above, it is enough to show this lower
%bound under the additionnal hypothesis that $\varphi \in {\mathcal H}_\theta$.
%We make this assumption in the sequel.
%$k_\varphi :=   \dot{\varphi}(s) - b(\varphi(s), \sigma_s)\in \mathcal H.$\\
 We set $r := \frac{1}{2} \inf_{s \in [0,1]} (\varphi(s) - 2 \sqrt{s})
>0.$

 From our assumptions on $\varphi,$ there exists $ \delta >0$ small enough such that:
\begin{equation} \label{P}
 \forall \chi \in B(\varphi, \delta), \forall \mu \in  \mathbb B_r(\sigma, \alpha), \forall s \in ]0,1] \,\textrm{ and } \,y \in \textrm{supp}(\mu_s), \chi(s) - y \geq \frac{r}{4}. 
\end{equation}

 For $h\in \mathcal H$ and $(\varphi, \mu) \in
 C_\T([0,1]; \R) \times C([0,1]; \mathcal P(\R)) $ such that for all $t
 \in [0,1],$  $\varphi(t) > r(\mu_t),$ we can define
  \begin{eqnarray}
  G_N(\varphi, \mu; h) &=& h(1) \varphi(1)  - h(0) \varphi(0) -
   \int_0^1 \varphi(s) \dot{h}(s) ds - \int_0^1 b_N(\varphi(s), \mu_s)
   h(s) ds,  \label{defG} \nonumber\\
   F_N(\varphi, \mu; h)  &:=& G_N(\varphi, \mu; h) - \frac{1}{2} \int_0^1 h^2(s) ds \label{defiF}\
  \end{eqnarray}
where we recall that $\displaystyle b_N = \frac{N-1}{N} b.$ 
  
 \noindent
  Therefore, from \eqref{defimart} and \eqref{EDSvp},  we have 
  \begin{equation} \label{martingale}
 M_1^h = \exp(NF_N(\lambda_1, \nu_N; h)).
 \end{equation}

\noindent
%$\nu_N = \frac{1}{N-1} \sum_{i=2}^N \delta_{\lambda_i}$.
We get
\begin{eqnarray*}
\PP(\lambda_1 \in B(\varphi, \delta)) & \geq & \PP\left(\lambda_1 \in
B(\varphi, \delta); \nu_N \in \mathbb B_r(\sigma, \alpha)\right) \\
&=& \E\left({\bf 1}_{\lambda_1 \in B(\varphi, \delta);  \nu_N \in
  \mathbb B_r(\sigma, \alpha)} \frac{M^{ k_\varphi}_1}{M^{ k_\varphi}_1}\right) \\
&=& \E^{ k_\varphi}\left({\bf 1}_{\lambda_1 \in B(\varphi, \delta);
  \nu_N \in \mathbb B_r(\sigma, \alpha)} \exp(-N F_N(\lambda_1, \nu_N; k_\varphi))\right) \\
&\geq &   \exp\left(-N \sup_{(\psi, \mu) \in C_{\alpha, \delta,r}}
F_N(\psi, \mu; k_\varphi)\right)\\
&& \times \,\PP^{ k_\varphi}\left(\lambda_1 \in B(\varphi,
\delta);  \nu_N \in \mathbb B_r(\sigma, \alpha)\right)
\end{eqnarray*}
where \begin{equation} \label{defC}
C_{\alpha, \delta, r} =  B(\varphi, \delta)
\times \mathbb B_r(\sigma, \alpha). \end{equation}

\noindent
Therefore
\begin{eqnarray}
\lefteqn{ \liminf_{N \rightarrow \infty} \frac{1}{N}  \ln\PP(\lambda_1 \in B(\varphi, \delta)) \geq - \sup_{(\psi, \mu) \in C_{\alpha, \delta,r}}  F(\psi, \mu; k_\varphi)  } \nonumber \\
&&  + \qquad \liminf_{N \rightarrow \infty} \frac{1}{N} \ln\PP^{k_\varphi}(\lambda_1 \in
B(\varphi, \delta);  \nu_N \in \mathbb B_r(\sigma, \alpha)). \label{borneinf}
\end{eqnarray}
 From the property \eqref{P} above, the fonction $(\psi, \mu) \mapsto
 F(\psi, \mu;k_\varphi)$ is continuous on $C_{\alpha, \delta,r}$ and we
 checked in Lemma \ref{IetJ} that $F(\varphi, \sigma; k_\varphi) = I_\T(\varphi).$ 

Moreover, from Proposition \ref{propPh}, we get that the last term in \eqref{borneinf} is equal to zero.

\noindent
We have thus obtained that for $\varphi \in {\mathcal H}_\T$, 
\begin{equation} \label{borneinf1}
 \lim_{\delta \downarrow 0} \liminf_{N \rightarrow \infty} \frac{1}{N}  \ln\PP(\lambda_1^{\T,N} \in B(\varphi, \delta)) \geq - I_\T(\varphi).
  \end{equation}
%It remains  to extend the above inequality to all $\varphi$ such that
%$I_\T(\varphi) < \infty$. This follows directly from the Lemma \ref{lemapprox}.

\subsection{Lower bound for a null initial condition}

We want to show Proposition \ref{lowerbound} - Eq.\eqref{eq-lowerbound}  under the assumption that $\theta =0.$

\smallskip\noindent
%Again from Lemma \ref{lemapprox}, we can assume that 
Let $\varphi \in {\mathcal H}_0$ and  we set $t_0(\varphi)$ as defined in \eqref{defH}. 
Then, $k_\varphi = 0$ on $[0, t_0(\varphi)]$.
We choose $\varepsilon$ given by Lemma \ref{lempos} and
 we denote by $r := r(\varepsilon) =  \frac{1}{2} \inf_{s \in [t_0(\varphi)+ \varepsilon,1]} (\varphi(s) - 2 \sqrt{s}) >0$. 
As in the case when $\theta >0$, we perform a change of measure via the martingale $M^{k_\varphi}$.
 Recall that $F_N$ is defined by \eqref{defiF}. We define $F_N^{(\varepsilon)}$ by
 \begin{eqnarray*}
 F_N^{(\varepsilon)}(\varphi, \mu; k_\varphi) &= &k_\varphi (1) \varphi(1)  - \int_{t_0(\varphi)}^1 \varphi(s) \dot{k}_\varphi(s) ds - \int_{t_0(\varphi) + \varepsilon}^1 b_N(\varphi(s), \mu_s) k_\varphi(s) ds \\
 && \qquad  \qquad \qquad - \frac{1}{2} \int_0^1 k_\varphi ^2(s) ds,
 \end{eqnarray*}
 in other words,
 $$ F_N(\varphi, \mu; k_\varphi )  = F_N^{(\varepsilon)}(\varphi, \mu; k_\varphi )  - \int_{t_0(\varphi)}^{t_0(\varphi) +\varepsilon} b_N(\varphi(s), \mu_s) k_\varphi (s) ds.$$

  Therefore, for such $\varepsilon$, $F_N \leq F_N^{(\varepsilon)}$ and we obtain (as in the previous subsection)
\begin{multline*}
   \PP(\lambda_1 \in B(\varphi, \delta))  \geq   \exp\left(-N \sup_{(\psi,
    \mu) \in C_{\alpha, \delta,r}} F_N^{(\varepsilon)}(\psi, \mu; k_\varphi )\right)\\
\PP^{k_\varphi}(\lambda_1 \in B(\varphi, \delta);  \nu_N \in \mathbb B_r(\sigma, \alpha)) 
\end{multline*}
    where $C_{\alpha, \delta, r}$ is defined in \eqref{defC}
and, using Proposition \ref{propPh}, 
$$ \liminf_{N \rightarrow \infty}\frac{1}{N} \ln \PP(\lambda_1 \in B(\varphi, \delta)) \geq - \sup_{(\psi, \mu) \in  C_{\alpha, \delta,r}} F_N^{(\varepsilon)}(\psi, \mu; k_\varphi ).$$
Now, for $\delta= \delta(\varepsilon)$ small enough, $F_N^{(\varepsilon)}$ is continuous on $ C_{\alpha, \delta,r}$, since
$$ \forall \psi \in B(\varphi, \delta), \forall \mu \in  \mathbb B_r(\sigma, \alpha), \forall s \in [ t_0(\varphi) + \varepsilon ,1]
\,\textrm{ and }\, y \in \textrm{supp}(\mu_s), \psi(s) - y \geq \frac{r}{4}
$$
 and therefore
$$\lim_{\delta \rightarrow 0} \liminf_{N \rightarrow \infty} \frac{1}{N} \ln \PP(\lambda_1 \in B(\varphi, \delta))\geq - F^{(\varepsilon)}(\varphi, \sigma; k_\varphi ),$$
where
\begin{multline*}
F^{(\varepsilon)}(\varphi, \sigma; k_\varphi ) = \frac{1}{2} \int_{t_0(\varphi) + \varepsilon}^1 (\dot{\varphi}(s) - b(\varphi(s), \sigma_s))^2 ds  - \int_{t_0(\varphi)}^{t_0(\varphi) +\varepsilon} \varphi(s) \dot{k}_\varphi (s) ds \\
- \frac{1}{2} \int_0^\varepsilon k_\varphi ^2(s) ds. 
\end{multline*}

This last quantity tends to $I_0(\varphi)$ as $\varepsilon$ tends to 0. \hfill$\Box$
%%%%%%%%%%%%%% a ajouter version Mylene
%%%%%%% THE UPPER BOUND

 \section{The upper bound}
\label{sec:upper}

We first prove the following
\begin{proposition}
\label{wub1}
Let $\T \geq 0$ and $\varphi \in C_\T([0,1]; \R)$ such that there exists $t_0 \in [0,1]$ so
that
$\varphi(t_0) < 2\sqrt {t_0}.$ Then
$$ \lim_{\delta \downarrow 0 } \lim_{N \rightarrow \infty}
\frac{1}{N} \ln \PP(\lambda_1^{\T,N} \in B(\varphi, \delta) ) = - \infty.$$ 
\end{proposition}

We proceed as in \cite{BDG}. From \cite{CDG}, we know that the
process $\mu_N$ satisfies a LDP in the scale $N^2$ with a good
rate function whose unique minimizer is the semicircular process
$\sigma$
for which we know that the support of $\sigma_t$ is $[-2\sqrt t, 2\sqrt t].$

Let $\delta_0 = 2\sqrt {t_0}- \varphi(t_0).$ By continuity of $\varphi, $
there exists $\varepsilon >0,$ such that for any $t \in [t_0-\varepsilon,
  t_0+\varepsilon],$
$\varphi(t) <2\sqrt t - \frac{\delta_0}{2}.$

\smallskip\noindent
 For any  $t \in [t_0-\varepsilon,
  t_0+\varepsilon],$ there exists $f_t$ such that $f_t(y) = 0$ if $y \leq
\varphi(t)$ and $\int f_t(x) d\sigma_t(x) >0.$
We let $F :=  \{\mu \in  C([0,1]; \mathcal P(\R))/\int f_t(x) d\mu_t(x)=0 \,\,\forall t \in [t_0-\varepsilon,
  t_0+\varepsilon]  \},$ which is a closed set.

\smallskip\noindent
 For any $\delta < \frac{\delta_0}{2},$ 
$$ \PP(\lambda_1 \in B(\varphi, \delta)) \leq \PP( \mu_N \in F).$$
As $\sigma \notin F, $
$\limsup_{N \rightarrow \infty} \frac{1}{N^2} \ln    \PP( \mu_N \in
F) < 0,$ which gives the Proposition. \hfill $\Box$

\smallskip\noindent 
We thus consider the case where $\varphi(t) \geq 2 \sqrt{t}$ and as a first step, we prove the upper bound for a function $\varphi$ which satisfies $\varphi(t) > 2 \sqrt{t}$
 for all $t \in [0,1]$ (this implies in particular that $\T >0$).
 
\subsection{The upper bound for functions $\varphi$ well separated from $t \mapsto 2\sqrt t$}
\begin{proposition}
\label{wub2}
Let   $\varphi \in C_\theta([0,1]; \R)$ such that 
%$\varphi$ is absolutely continuous and 
 for any $t \in [0,1],$ $\varphi(t) > 2\sqrt t.$
 Then
\begin{equation} \label{wub3}
\lim_{\delta \downarrow 0 } \limsup_{N \rightarrow \infty}
\frac{1}{N} \ln \PP(\lambda_1^{\T,N} \in B(\varphi, \delta) ) \leq - I_\T(\varphi).
\end{equation}
\end{proposition}
%Before proving this theorem, we give some preliminary results. \\
%Let $\varphi \in C_\T([0,1]; \mathbb{R})$ such that  $\varphi(t) > 2 \sqrt{t}$ for all $t \in [0,1].$  
%We recall that  $r := \frac 1 2 \inf_{t} (\varphi(t) - 2 \sqrt{t}) >0$. 
The general idea to prove \eqref{wub3} is, as for the lower bound, to introduce the exponential martingales $M^h_t$ (see \eqref{defimart} and \eqref{martingale}) and to optimize on $h$.
\begin{equation} \label{essai}
\PP(\lambda_1 \in B(\varphi, \delta))  \leq  \PP\left(\lambda_1 \in
B(\varphi, \delta); \nu_N \in \mathbb B(\sigma, \alpha)\right) + \PP( \nu_N \not\in \mathbb B(\sigma, \alpha))
\end{equation}
and the second term is exponentially negligible in scale $N$ since $\nu_N$ satisfies  a large deviation principle in scale $N^2$ (this would be not true if we replace $\mathbb B(\sigma, \alpha)$ by $\mathbb B_r(\sigma, \alpha)$). Now,
\begin{eqnarray*}
\PP\left(\lambda_1 \in
B(\varphi, \delta); \nu_N \in \mathbb B(\sigma, \alpha)\right)&=& \E\left({\bf 1}_{\lambda_1 \in B(\varphi, \delta);  \nu_N \in
  \mathbb B(\sigma, \alpha)} \frac{M^h_1}{M^h_1}\right) \\
&\leq &   \exp\left(-N \inf_{(\psi, \mu) \in B(\varphi, \delta)\times \mathbb B(\sigma, \alpha) }
F_N(\psi, \mu;h)\right).
\end{eqnarray*}
Unfortunately, the fonction $F$ (or $F_N$) is not continuous on $ B(\varphi, \delta)\times B(\sigma, \alpha)$ and we cannot conclude that 
$$ \lim_{\delta \rightarrow 0, \alpha \rightarrow 0}  \inf_{(\psi, \mu) \in B(\varphi, \delta)\times \mathbb B(\sigma, \alpha) }
F_N(\psi, \mu; h) = F(\varphi, \sigma; h).$$
Therefore, the strategy of the proof is first to prove that with high probability,  only a finite number (say $K$) of eigenvalues can deviate strictly above $t \mapsto 2\sqrt t$ and then to introduce exponential martingales depending on $\lambda_1, \ldots, \lambda_K$  involving a functional
$F(\lambda_1, \ldots, \lambda_K,\mu_N^{(1)}, \ldots, \mu_N^{(K)})$  (see \eqref{martingale1}) which is now continuous on the sets that we are considering. 

We first prove 
\begin{proposition} \label{propK}
 For any $\eta >0,$ $L>0,$ there exists
%Let $\varphi$ such that $I_\T(\varphi) < \infty$ and let $\eta$ as above. There exists 
$K:= K(\eta,L)$ (independent of $N$) such that
\begin{equation}
\limsup_{N \rightarrow \infty} \frac{1}{N} \ln\PP( \exists t \in [0,1], \lambda_{K+1}^{\T,N}(t) > 2 \sqrt{t} + \eta) \leq -L.
\end{equation}
\end{proposition}
We will first prove 
a fixed time version of the same result stated in the following lemma:
\begin{lemma}
 \label{fixtime}
%Let $\varphi$ such that $I_\T(\varphi) < \infty$ and let $\eta$ as above.
%If we denote by $J$ the rate function for $\lambda_1^{N,0}(1)$ and let $K(\eta,\varphi) := \lfloor \frac{2 I_\T(\varphi)}{J(2+\eta/2)}\rfloor.$
 For any $\eta >0,$ $L>0,$ there exists
%Let $\varphi$ such that $I_\T(\varphi) < \infty$ and let $\eta$ as above. There exists 
$K:= K(\eta,L)$  such that for any  $t \in [0,1],$ 

$$ \limsup_{N \rightarrow \infty} \frac{1}{N} \ln\PP( \lambda_{K+1}^{\theta,N}(t) > 2\sqrt{t}  + \eta) \leq -L.$$
\end{lemma}

\noindent{\bf Proof:}
The first observation is that, as $H_N^\T(t) = H_N^0(t) + \textrm{diag}(\T, 0, \ldots,0),$ with $\T \ge 0,$ by Weyl's inequalities,
$ \lambda_{K+1}^{\T,N}(t) \le \lambda_{K}^{0,N}(t)$ so that
\begin{multline*}
 \PP ( \lambda_{K+1}^{\T,N}(t) > 2\sqrt{t}+\eta) \le \PP \left( \lambda_{K}^{0,N}(t) > 2\sqrt{t}+\eta\right) = 
\PP\left( \lambda_{K}^{0,N}(1) > 2+\frac{\eta}{\sqrt{t}}\right) \\
\leq \PP ( \lambda_{K}^{0,N}(1) > 2+\eta).
\end{multline*}
Therefore, Lemma \ref{fixtime} will be a direct consequence of the fact that for any $p \ge 1,$
the law of $(\lambda_{1}^{0,N}(1), \ldots, \lambda_{p}^{0,N}(1))$ satisfies a LDP in the scale $N$ with good rate function
$$ F: (x_1, \ldots, x_p) \mapsto 1_{x_1 \ge x_2 \ge \ldots \ge x_p} \sum_{i=1}^p  K_0(x_j),$$
with $K_0$ the individual rate function at time 1 as defined in Theorem \ref{pgdflou}.
This is a particular case of Theorem 2.10 in \cite{BGM} in the case when the potential $V$ is just Gaussian
($V(x)=x^2$) therein.

\smallskip\noindent
 From this, if we define $K \ge \frac{L}{K_0(2+\eta)}$, we deduce that, 
\begin{multline*}
 \limsup_N \frac{1}{N} \ln\PP( \lambda_{K+1}^{\T,N}(t) > 2\sqrt{t}  + \eta) 
\le \limsup_N \frac{1}{N} \ln\PP( \lambda_{K}^{0,N}(1) > 2  + \eta) \\
=  \limsup_N \frac{1}{N} \ln\PP( \lambda_{1}^{0,N}(1) > 2  + \eta,\ldots,\lambda_{K}^{0,N}(1) > 2  + \eta)\\
\hspace{8cm} \le - K K_0(2+\eta) \le -L. \hfill\Box\\
\end{multline*}

\noindent
{\bf Proof of Proposition \ref{propK}:} %Let  $0<c<1/2$ and a subdivision $(t_k^{N}, k \leq N^c)$ of $[0,1]$ with $t_{k+1}^{(N)} - t_k^{(N)} = N^{-c}.$ 
%Let $K:= K(\eta,C)$ defined in Lemma \ref{fixtime} above, then 
We fix $\eta >0$ and $L>0.$ Let $R$ be such that $\frac{1}{10}\eta^2R >18L$ and choose
a subdivision $(t_k)_{1 \le k \le R}$ of the interval $[0,1]$ such that for all $1 \le k \le R,$ $|t_k-t_{k+1}| \le \frac{2}{R}.$ 
Now, for any $K \in \mathbb N^*$ 
\begin{eqnarray}
\PP( \exists t \in [0,1], \lambda_{K+1}(t) &> &2 \sqrt{t} + \eta) =  \PP[ \cup_k (\exists t \in [t_k,t_{k+1}], \lambda_{K+1}(t) > 2 \sqrt{t} + \eta) ] \nonumber\\
&\leq & R \max_{1\le k \le R} \PP(  \exists t \in [t_k,t_{k+1}], \lambda_{K+1}(t) > 2 \sqrt{t} + \eta). \label{eqR}
\end{eqnarray}
Then
\begin{eqnarray*}
\lefteqn{
\PP(  \exists t \in [t_k,t_{k+1}], \lambda_{K+1}(t) > 2 \sqrt{t} + \eta) 
\leq \PP\left(  \lambda_{K+1}(t_k)> 2 \sqrt{t_k} + \frac{\eta}{2}\right)\qquad \qquad \qquad \qquad  }  \\
\qquad &&   \qquad    + \PP\left(  \exists t \in [t_k,t_{k+1}], \lambda_{K+1}(t) > 2 \sqrt{t} + \eta \ ;  \lambda_{K+1}(t_k) \leq  2 \sqrt{t_k} + \frac{\eta}{2}\right) \\
 \qquad& \leq &\PP\left(  \lambda_{K+1}(t_k) > 2 \sqrt{t_k} + \frac{\eta}{2}\right) + \PP\left( \sup_{t_k \leq t < t_{k+1}} \vert \lambda_{K+1}(t) - \lambda_{K+1}(t_k)\vert \geq \frac{\eta}{3} \right).
\end{eqnarray*}
\smallskip

 From Lemma \ref{fixtime}, we can find $K:= K(\eta, L)$ such that
$$\limsup_N \frac{1}{N} \ln\PP\left(  \lambda_{K+1}(t_k) > 2 \sqrt{t_k} + \frac{\eta}{2}\right)\leq- L.$$
 From Lemma \ref{lemexptight} applied for $p=K(\eta, L),$ 
$$\limsup_N \frac{1}{N} \ln\PP\left( \sup_{t_k \leq t < t_{k+1}} \vert \lambda_{K+1}(t) - \lambda_{K+1}(t_k)\vert 
\geq \frac{\eta}{3} \right)=-L.$$
As $R$ is independent of $N,$  \eqref{eqR} gives the lemma. 
\hfill $\Box$

%%%%%%%%%%%%%%%%%%%%%%%%%
\noindent

%%%%
%We end this preliminary part by a lemma:
\begin{lemma} \label{mesemp}
Let $K$ be fixed as in Proposition \ref{propK}. Let $j \in \{1, \ldots K\}$. We denote by $\mu_N^{(j)}$ the spectral measure of the $N-j$ smallest eigenvalues
$\displaystyle\mu_N^{(j)} = \frac{1}{N-j} \sum_{p=j+1}^N \delta_{\lambda_p}.$
Then, 
$$
d(\mu_N^{(j)} , \mu_N) \leq \frac{2K}{N}.$$
Therefore, if $\mu_N \in \mathbb{B}(\sigma, \alpha)$, $\mu_N^{(j)} \in \mathbb{B}(\sigma, 2 \alpha)$ for $N \geq N_0$.
\end{lemma}

\noindent{\bf Proof:} Let $f \in {\mathcal F}_{Lip}$,
$$ \mu_N^{(j)} (f) - \mu_N(f) = \frac{j}{N(N-j)} \sum_{p>j} f(\lambda_p) - \frac{1}{N} \sum_{p=1}^j f(\lambda_p)$$
and 
$$ \vert  \mu_N^{(j)} (f) - \mu_N(f) \vert \leq  \frac{j}{N} +  \frac{j}{N} \leq  \frac{2K}{N}. $$
\hfill $\Box$
%%%%%%%%%%%%%

%%%%%%%%%%%%%%%%
\noindent
{\bf Proof of Proposition \ref{wub2}:} 
Let   $\varphi \in C_\theta([0,1]; \R)$ such that 
%$\varphi$ is absolutely continuous and 
 for any $t \in [0,1],$ $\varphi(t) > 2\sqrt t.$ We recall that $r =
 \frac{1}{2} \inf(\varphi(t) - 2\sqrt t)$\smallskip

%a) We consider $\varphi$ such that $I_\T(\varphi) < \infty$. \\
 For $K \in \N^*,$ let $\delta >0$ such that $\delta < \frac{r}{4K}$ and $\alpha >0$. We have
\begin{equation} \label{3term}
 \PP( \lambda_1 \in B(\varphi, \delta)) \leq 
\PP(A_{N, \delta, \alpha,K} ) 
+ \PP(\exists t \in [0,1], \lambda_{K+1}(t) > 2 \sqrt{t} + r) + \PP( \mu_N \not\in \mathbb B(\sigma, \alpha)) 
\end{equation}
with 
$$ A_{N, \delta, \alpha,K}  := \left\{\lambda_1 \in B(\varphi, \delta); \forall p>K, \forall t, \lambda_p(t) \leq 2\sqrt{t} + r; \ \mu_N \in \mathbb B(\sigma, \alpha)\right\}.$$

%From the previous results,
%$$\limsup_N \frac{1}{N} \ln(\PP( \exists t \in [0,1], \lambda_{K+1}(t) > 2 \sqrt{t} + \eta)) \leq -2 I_\T(\varphi)$$
%and
%$$\limsup_N \frac{1}{N} \ln(\PP( \mu_N \not\in \mathbb B(\sigma, \alpha))) = - \infty.$$
%We denote by $A_{N, \delta, \alpha}$ the event in \eqref{3term}.
 
 For each $t \in [0,1]$, if we are on $A_{N, \delta, \alpha,K}$, there exists at least a gap between two eigenvalues (among the $K$ largest ones) larger than $\delta$.
%For $1 \le k \le R,$ we can define a random variable $I_k$ with values in $\{1, \ldots, K\}$ by
%$$I_k := \inf \{ i \leq K, \lambda_{i}(t_k) - \lambda_{i+1}(t_k)> \delta \}.$$
To proceed, we will choose a subdivision $(t_k)_{1 \leq k \leq R}$ as in the proof of Proposition \ref{propK} and we will decompose
the event $A_{N, \delta, \alpha,K}$ according to the location of the gap of size $\delta.$ 
More precisely, we have
\begin{equation}
 \label{decompo1}
A_{N, \delta, \alpha,K} \subset \bigcup_{{\bf i} \in\{1, \ldots K\}^R} A_{N, {\bf i}, \delta, \alpha}
\end{equation}
where, for ${\bf i} = (i_1,\ldots,i_R)\in \{1, \ldots,K\}^R,$
\begin{multline*}
 A_{N, {\bf i}, \delta, \alpha} = \left\{ \lambda_1 \in B(\varphi, \delta), \forall k\leq R, \forall i< i_k,  \lambda_{i}(t_k) - \lambda_{i+1}(t_k)\leq \delta, \right.\\ \left.\lambda_{i_k}(t_k) - \lambda_{i_k+1}(t_k)> \delta, \mu_N \in \mathbb B(\sigma, \alpha)  \right\}.
\end{multline*}
As, for $i \leq i_k$, $ \lambda_i(t_k) \geq \varphi(t_k) - i\delta$, 
\begin{multline*}
 A_{N, {\bf i}, \delta, \alpha} \subset \left\{ \lambda_1 \in B(\varphi, \delta), \forall k\leq R, \forall i< i_k,  \lambda_{i}(t_k) \geq \varphi(t_k)  -i  \delta,\right.\\ \left. \lambda_{i_k}(t_k) - \lambda_{i_k+1}(t_k)> \delta, \mu_N \in \mathbb B(\sigma, \alpha)   \right\}.
\end{multline*}
%Up to some error term exponentially negligible in scale $N$ (see Lemma \ref{lemexptight} applied for $\delta = N^{-c}$),
%we can replace this event  by the event
Now, we choose $R$ such that $\sup_{|t-s| \leq \frac{2}{R}} |\varphi(s) - \varphi(t)| \leq \frac{\delta}{6}$
and the subdivision such that $|t_k - t_{k+1}|\leq \frac{2}{R}.$ If we let 
\begin{multline*}
  B_{N, {\bf i}, \delta, \alpha} = \left\{\lambda_1 \in B(\varphi, \delta), \forall i< i_k,  \forall t \in [t_k, t_{k+1}[,
 \lambda_{i}(t) \geq \varphi(t) - \left(i+ \frac{1}{3}\right)\delta, \right. \\
\left. \lambda_{i_k}(t) - \lambda_{i_k+1}(t)>\frac{2}{3} \delta,
\mu_N \in \mathbb B(\sigma, \alpha)  \right\},
\end{multline*}
then 
\begin{equation}
  \label{decompo2}
 A_{N, {\bf i}, \delta, \alpha}  \subset B_{N, {\bf i}, \delta, \alpha} \bigcup
\left\{\exists k, \exists t \in [t_k,t_{k+1}], \exists i \leq i_k +1,  | \lambda_{i}(t) - \lambda_{i}(t_k)|> \frac{\delta}{6}  \right\}.
\end{equation}
The second term will again be controlled by Lemma \ref{lemexptight} and 
we now need to work on  $B_{N, {\bf i}, \delta, \alpha}.$ Our goal is to show that  for any $K \in \N^*,$  any $h \in \mathcal H,$ and any subdivision $(t_k)_{1 \le k \le R}$ of $[0,1],$
\begin{equation}\label{limBN}
 \lim_{\delta \rightarrow 0}  \lim_{\alpha \rightarrow 0} \limsup \frac{1}{N} \ln \PP(B_{N, {\bf i}, \delta, \alpha} ) \leq - F(\varphi, \sigma; h).
\end{equation}
The idea is, as for the lower bound to make a change of measure given by a martingale but this time
it will depend not only on $\lambda_1$ but on all the eigenvalues near $\varphi$ above the gap of size $\delta.$
Then the average of these eigenvalues will be near $\varphi$ and the variance of their Brownian part will be smaller 
than for an individual eigenvalue.

More precisely, let $j \leq K$ and $X_j(t) := \frac{1}{j} \sum_{i=1}^j \lambda_i(t)$ is a solution of the SDE
\begin{equation}
dX_j(t) = \frac{1}{\sqrt{N}} \frac{1}{j} \sum_{i\leq j} dB_i(s) + \frac{1}{N}\frac{1}{j} \sum_{i\leq j} \sum_{p>j} \frac{1}{\lambda_i(t) - \lambda_p(t)} dt.
 \end{equation}
We denote by $dB^{j}(s) = \frac{1}{\sqrt{j}} \sum_{i=1}^j dB_i(s),$ which is  a standard Brownian motion. 
Let $h \in \mathcal H,$ we define
 the exponential martingale 
\begin{multline*}
 \widetilde  M_t^h = \exp \left[ N \left( \int_0^t \sum_{k\le R} h(s) 1_{[t_k, t_{k+1}[}(s)  \frac{1}{\sqrt{N}} \frac{1}{\sqrt{i_k}} 
  dB^{i_k} (s) \right.\right.\\\left.\left.- \frac{1}{2} \int_0^t \sum_{k \le R}  1_{[t_k, t_{k+1}[}(s) \frac{1}{i_k} h^2(s) ds\right) \right].
\end{multline*}
\begin{eqnarray*}
\frac{1}{\sqrt{N}} \frac{1}{\sqrt{i_k}} dB^{i_k} (t) &=& dX_{i_k}(t) - \frac{1}{N}\frac{1}{i_k} \sum_{i\leq i_k} \sum_{p>i_k} \frac{1}{\lambda_i(t) - \lambda_p(t)} dt \\
&=& dX_{i_k}(t) - \frac{N-i_k}{N}\frac{1}{i_k} \sum_{i\leq i_k}\int \frac{(\mu_N^{(i_k)})_t(dx)}{\lambda_i(t) - x} \ dt 
\end{eqnarray*}
from the definition of the measures $\mu_N^{(i)}$ given in Lemma \ref{mesemp}.  Thus
\begin{eqnarray*}
\widetilde M_1^h &= &\exp \left[ N \sum_k \left( [h_s X_{i_k}(s)]_{t_k}^{t_{k+1}} - \frac{N-i_k}{N}\frac{1}{i_k} \sum_{i\leq i_k}\int_{t_k}^{t_{k+1}} \int \frac{\mu_t^{(i_k)}(dx)}{\lambda_i(t) - x} h(t) \ dt \right. \right. \\
&& -
 \int_{t_k}^{t_{k+1}} \dot{h}(s) X_{i_k}(s) ds  \left. \left. - \frac{1}{2i_k} \int_{t_k}^{t_{k+1}}h^2(s) ds\right) \right].
 \end{eqnarray*}
 We recall from Lemma \ref{mesemp} that if $\mu_N \in \mathbb B(\sigma, \alpha)$, $\mu_N^{(i_k)} \in \mathbb B(\sigma, 2\alpha)$. \\
 $\widetilde M^h_1$ can be written as a functional
 \begin{equation} \label{martingale1}
\widetilde M_1^h = \exp\left( NF_N(\lambda_1, \ldots \lambda_K, \mu_N^{(1)}, \ldots, \mu_N^{(K)}; h)\right).
 \end{equation}
 We denote by 
 \begin{eqnarray*}
 \Lambda_{\bf{i}, \delta, \alpha} &= &\{(\psi_1, \ldots \psi_K, \nu_1,\ldots  \nu_K): \psi_1 \in B(\varphi, \delta),\\
&& \forall k\le R,
 \forall i< i_k,  \forall t \in [t_k, t_{k+1}[,\psi_{i}(t) \geq \varphi(t) - \left(i+ \frac 1 3\right)\delta, \\
 && \;  \psi_{i_k}(t) - \psi_{i_{k+1}}(t)>\frac{2}{3} \delta; \nu_i  \in \mathbb B(\sigma, 2\alpha), \textrm{supp}(\nu_{i_k}(.))  \subset ]-\infty, \psi_{i_k+1}(.) \}
 \end{eqnarray*}
where in the above set, the functions are such that $\psi_1 \geq \psi_2 \ldots \geq \psi_K$.
We denote by $\underline{\psi} = (\psi_1, \ldots \psi_K)$ and $\underline{\nu} = (\nu_1, \ldots, \nu_K)$.
Then,
\begin{eqnarray*}
\PP( B_{N, {\bf i}, \delta, \alpha}) &= &\E\left[ 1_{ B_{N, {\bf i}, \delta, \alpha}} \frac{M_1^h}{M_1^h}\right] \\
&\leq & \exp\left(- N \inf_{ (\underline{\psi} , \underline{\nu} ) \in \Lambda_{\bf{i}, \delta, \alpha} } F_N(
\underline{\psi} , \underline{\nu} ; h\right) \E[M_1^h] \\
&\leq & \exp\left(-N \inf_{ (\underline{\psi} , \underline{\nu} ) \in \Lambda_{\bf{i}, \delta, \alpha} } F_N(
\underline{\psi} , \underline{\nu} ; h)\right)
\end{eqnarray*}
and
$$\frac{1}{N} \ln\PP( B_{N, {\bf i}, \delta, \alpha}) \leq - \inf_{ (\underline{\psi} , \underline{\nu} ) \in \Lambda_{\bf{i}, \delta, \alpha} } F_N(\underline{\psi} , \underline{\nu} ; h),$$
$$\limsup \frac{1}{N} \ln\PP( B_{N, {\bf i}, \delta, \alpha}) \leq - \inf_{ (\underline{\psi} , \underline{\nu} ) \in \Lambda_{\bf{i}, \delta, \alpha} } F_{\bf i}(\underline{\psi} , \underline{\nu} ; h)$$
where
\begin{eqnarray*}
F_{\bf i}(\underline{\psi} , \underline{\nu} ; h)  &=& \sum_k\left[ [h_s \Psi_{i_k}(s)]_{t_k}^{t_{k+1}} -\frac{1}{i_k} \sum_{i\leq i_k}\int_{t_k}^{t_{k+1}} \int \frac{(\nu_{i_k})_t(dx)}{\psi_i(t) - x} h(t) \ dt \right.\\
&& -
 \int_{t_k}^{t_{k+1}} \dot{h}(s) \Psi_{i_k}(s) ds-\left. \frac{1}{2i_k} \int_{t_k}^{t_{k+1}}h^2(s) ds\right]
 \end{eqnarray*}
 with $\Psi_j = \frac{1}{j} \sum_{i\leq j} \psi_i$. 

 Let us take $\alpha \rightarrow 0$. The function $\underline{\nu} \mapsto F_{\bf i}(\underline{\psi} , \underline{\nu} ; h) $ is continuous on the set
 $ \Lambda_{\bf{i}, \delta, \alpha}$ since for $i \leq i_k$, 
 $$ \psi_i(t) - x \geq \frac{2}{3} \delta \quad  \forall x \in \textrm{supp}( (\nu_{i_k})_t).$$
 We obtain
 $$\lim_{\alpha \rightarrow 0} \limsup \frac{1}{N} \ln\PP( B_{N, {\bf i}, \delta, \alpha})\leq - \inf_{ \underline{\psi}  \in \Lambda_{\bf{i}, \delta } } F_{\bf i}(\underline{\psi} , \underline{\sigma} ; h)$$
 where $\underline{\sigma} = ( \sigma, \ldots \sigma)$ and $ \Lambda_{\bf{i}, \delta}$ is defined as in $ \Lambda_{\bf{i}, \delta, \alpha}$ without the conditions on $\nu_i$.
 Now, take $\delta \rightarrow 0$, the above functional is continuous in $\underline{\psi}$ and 
 \begin{equation}
\lim_{\delta \rightarrow 0}  \lim_{\alpha \rightarrow 0} \limsup \frac{1}{N} \ln\PP( B_{N, {\bf i}, \delta, \alpha}) \leq -  F_{\bf i}(\underline{\varphi} , \underline{\sigma} ; h)
\end{equation}
where
\begin{multline*}
 F_{\bf i}(\underline{\varphi} , \underline{\sigma} ; h) = h(1)\varphi(1) - h(0) \varphi(0) - \int_0^1 \int_\R \frac{\sigma_t(dx)}{\varphi(t) - x}h(t) dt - \int_0^t \dot h(s)\varphi(s)ds \\- \sum_k \frac{1}{2i_k} \int_{t_k}^{t_{k+1}}h^2(s) ds
\end{multline*}
and 
$$ - F_{\bf i}(\underline{\varphi} , \underline{\sigma} ; h)  \leq - F(\varphi, \sigma; h)$$ where 
$F$ is defined by \eqref{defF}.

 We have proved \eqref{limBN}.

We now  go back to the decompositions \eqref{decompo1} and \eqref{decompo2}. Let us first treat the case when $I_\T(\varphi) <\infty.$ We choose  $L= - 2I_\T(\varphi)$
and $K$ as given in Proposition \ref{propK} so that
$$\limsup_N \frac{1}{N} \ln\PP( \exists t \in [0,1], \lambda_{K+1}(t) > 2 \sqrt{t} + \eta) \leq -2 I_\T(\varphi).$$
Moreover
$$\limsup_N \frac{1}{N} \ln\PP( \mu_N \not\in \mathbb B(\sigma, \alpha)) = - \infty$$
and from Lemma \ref{lemexptight}, if we choose $R,$ the number of points of the subdivision such that $R > \frac{260 I_\T(\varphi)}{ \delta^2},$
$$ \limsup_N \frac{1}{N} \ln\PP\left( \exists k, \exists t \in [t_k,t_{k+1}], \exists i \leq i_k +1 | \lambda_{i_k}(t) - \lambda_{i_k}(t_k)|> \frac{\delta}{6} \right) 
\leq -2 I_\T(\varphi).$$

\noindent
We thus obtain, for any $h \in  \mathcal H,$
$$ \lim_{\delta \rightarrow 0} \limsup \frac{1}{N} \ln \PP(\lambda_1 \in B(\varphi, \delta) ) \leq - \inf( F(\varphi , \sigma ; h), 2 I_\T(\varphi)) .$$
Optimizing in $h$ gives
$$ \lim_{\delta \rightarrow 0}  \limsup_N \frac{1}{N} \ln \PP(\lambda_1 \in B(\varphi, \delta) ) \leq - I_\T(\varphi). $$

\noindent
In the case where $I_\T(\varphi) = \infty$, for any $L$, we can associate $K$  as in Proposition \ref{propK} such that
\begin{equation}
\limsup_N \frac{1}{N} \ln\PP( \exists t \in [0,1], \lambda_{K+1}(t) > 2 \sqrt{t} + \eta) \leq -L.
\end{equation}
In the same way as above, with $R > \frac{180L}{\delta^2},$ we then show that
$$ \lim_{\delta \rightarrow 0}  \limsup \frac{1}{N} \ln \PP(\lambda_1 \in B(\varphi, \delta) ) \leq -L $$
and since the left hand side does not depend on $L$,
$$ \lim_{\delta \rightarrow 0}  \limsup \frac{1}{N} \ln\PP(\lambda_1 \in B(\varphi, \delta) )= - \infty. $$
\hfill $\Box$

\noindent
 We now extend Proposition \ref{wub2} to any function $\varphi$ such that $\varphi(t) \geq 2 \sqrt{t}$. 
 %We assume that $\varphi(t) \not\equiv 2 \sqrt{t}$ since \eqref{wub3} is trivial for $\varphi= 2 \sqrt{.}$.
\subsection{The upper bound for functions $\varphi$ not well separated from $t \mapsto 2\sqrt t$}
\begin{proposition}
\label{wub4}
Let   $\varphi \in C_\theta([0,1]; \R)$ such that 
%$\varphi$ is absolutely continuous and 
 for any $t \in [0,1],$ $\varphi(t) \geq 2\sqrt t.$
 Then
$$
\lim_{\delta \downarrow 0 } \limsup_{N \rightarrow \infty}
\frac{1}{N} \ln \PP(\lambda_1^\T \in B(\varphi, \delta) ) \leq - I_\T(\varphi).$$
\end{proposition}

\noindent
{\bf Proof of Proposition \ref{wub4}:} For any $\epsilon >0,$
 let $J_\epsilon = \{ t \in [0,1], \varphi(t) > 2 \sqrt{t} + \epsilon \}.$ $J_\epsilon$ is an  open set in $[0,1]$ and $\overline J_\epsilon$ is compact so that
 we can find a set $V_\epsilon$ of the form $V_\epsilon = \cup_{i=1}^{N_\epsilon} ]a_i(\epsilon), b_i(\epsilon)[$ such that
 $$ \bar{J_\epsilon} \subset V_\epsilon \subset J_{\epsilon/2}.$$ 
% with $0 \leq a_1 < b_1 \leq a_2 < b_2 \ldots  \leq a_r < b_r \leq 1$. The intervals can be closed in $a_1 =0$ and $b_r = 1$.
% Let $0 < \epsilon < \frac{1}{4} \min(b_i -a_i)$ and consider the closed set $I_\epsilon := \cup_{i=1}^r [a_i+ \epsilon, b_i-\epsilon]$. 
Then, on $V_\epsilon$, $\varphi(t) > 2 \sqrt{t}$. For a function $f$ on $[0,1]$, we denote by $f\vert_A$ its restriction to a subset $A$ of $[0,1]$. Then,
 $$\PP(\lambda_1 \in B(\varphi, \delta)) \leq \PP(\lambda_1\vert_{V_\epsilon} \in B(\varphi\vert_{V_\epsilon} , \delta)).$$
 From Proposition \ref{wub2},
 $$ \lim_{\delta \downarrow 0 } \limsup_{N \rightarrow \infty}
\frac{1}{N} \ln \PP(\lambda_1\vert_{V_\epsilon} \in B(\varphi\vert_{V_\epsilon} , \delta) ) \leq - \sum_{i=1}^{N_\epsilon} I_\T(\varphi\vert_{[a_i(\epsilon), b_i( \epsilon)]} )$$
where 
$$ I_\T(\varphi\vert_{[a, b]}) = 
 \frac 1 2\int_a^b \left( \dot{\varphi}(s) -\frac{1}{2s}\left(\varphi(s) -
\sqrt{\varphi(s)^2-4s}\right)\right)^2 ds,$$
this quantity may be infinite.
Let $\epsilon \rightarrow 0$, by monotone convergence,
 $$ \lim_{\delta \downarrow 0 } \limsup_{N \rightarrow \infty}
\frac{1}{N} \ln \PP(\lambda_1  \in B(\varphi , \delta) ) \leq -   \frac 1 2\int_J \left( \dot{\varphi}(s) -\frac{1}{2s}\left(\varphi(s) -
\sqrt{\varphi(s)^2-4s}\right)\right)^2 ds$$
where $J= \{ t \in [0,1], \varphi(t) > 2 \sqrt{t} \}$ (the right-hand side can be infinite).

Assume that $\varphi$ is differentiable almost everywhere (a.e.)  on [$0,1].$
Since $\varphi(t) \geq 2 \sqrt{t}$, if $\varphi$ is differentiable in $s_0$ such that $\varphi(s_0) = 2 \sqrt{s_0}$: then, $\dot\varphi(s_0) = \frac{1}{\sqrt{s_0}}$ and
$$ \dot{\varphi}(s_0) -\frac{1}{2s_0}\left(\varphi(s_0) -
\sqrt{\varphi(s_0)^2-4s_0}\right) = 0.$$
Thus, 
$$ \int_J \left( \dot{\varphi}(s) -\frac{1}{2s}\left(\varphi(s) -
\sqrt{\varphi(s)^2-4s}\right)\right)^2 ds = I_\T (\varphi).$$
%%%%%%%
If $\varphi$ is not differentiable a.e., then $I_\T(\varphi) = \infty$. Consider first the case $\T >0$. From the lower semicontinuity of $I_\T$, for all $C >0$,  there exists $\epsilon$ such that 
$$B(\varphi, \epsilon) \subset \{ \psi; I_\T(\psi) > C\}.$$
Define $$ \left \{  \begin{array}{ll}\psi(t) = \varphi(t) & \rm{ on } \; \bar J_\epsilon \\
\psi(t) = 2 \sqrt{t} + \epsilon &  \rm{ on } \; (\bar J_\epsilon)^c  \end{array} \right. $$
Then, $\psi \in B(\varphi, \epsilon)$ and 
\begin{multline*}
\int_{(\bar J_\epsilon)^c} \left( \dot{\psi}(s) -\frac{1}{2s}\left(\psi(s) -
\sqrt{\psi(s)^2-4s}\right)\right)^2 ds \\
= \int_{(\bar J_\epsilon)^c} \left( \frac{1}{2s}\left(\epsilon -
\sqrt{\epsilon^2 + 4 \sqrt{s} \epsilon}\right)\right)^2 ds 
\leq  K\epsilon,
\end{multline*}
for some constant $K$. The last inequality follows from the fact that since $\T >0$, $(\bar J_\epsilon)^c \subset [a,1]$ for a strictly positive $a$.
Therefore, for $\epsilon$ small enough, $I_\T(\psi \vert_{\bar J_\epsilon}) = I_\T(\varphi  \vert_{\bar J_\epsilon})  \geq \frac{C}{2}.$
Moreover,  $I_\T(\varphi  \vert_{V_\epsilon}) \ge I_\T(\varphi  \vert_{\bar J_\epsilon})$ so that we get
$$\lim_{\delta \downarrow 0 } \limsup_{N \rightarrow \infty}
\frac{1}{N} \ln \PP(\lambda_1 \in B(\varphi, \delta)) \leq - \frac{C}{2}.$$
Since the inequality is true for all $C$,
$$\lim_{\delta \downarrow 0 } \limsup_{N \rightarrow \infty}
\frac{1}{N} \ln \PP(\lambda_1 \in B(\varphi, \delta))  = - \infty.$$
Now for $\T = 0$, if $I_0(\varphi \vert_{[a,1]}) < \infty$ for all $a>0$, then, $\varphi$ would be a.e. differentiable. 
Therefore, we can assume that there exists a $a$ such that $I_0(\varphi \vert_{[a,1]})  = \infty$ and argue as before, using that 
$\PP(\lambda_1 \in B(\varphi, \delta)) \leq \PP(\lambda_1\vert_{[a,1]} \in B(\varphi\vert_{[a,1]}, \delta)).$\hfill $\Box$

%%%%%%%%%%%%%%%%%
%In conclusion, \eqref{wub3} is satisfied for any $\varphi \in C_\T([0,1]; \R)$, proving the weak upper bound.
%Now, Theorem \ref{main} follows from the exponential tightness, the lower bound obtained in Proposition \ref{lowerbound} and the weak upper bound \eqref{wub3} (see \cite[Chapt. 4]{DZ}).

%%%%%%%%%%%%%%%%%%%%%%%%%%%%%%%%%%%%%%%%%%%%%%
%%%%%%%%%%%%%%%%%%%%%%%%%%%%%%%%%%%%%%%%%%%%%%%%%%%

\section{Contraction principle}
\label{sec:contract}

%In this section, we work under the assumption that Conjecture \ref{main} holds.
The goal of this section is to get from Theorem \ref{main}  a
new proof of Theorem \ref{pgdflou} concerning the deviations of the
largest eigenvalue at fixed time (say $t =1$).

\smallskip\noindent
Note that this fixed time result has been used in the preceding section for the proof of 
the upper bound in the case $\theta = 0,$ the goal here is to extend it to any $\theta >0.$\smallskip

\noindent
{\bf Proof of Theorem \ref{pgdflou} :}
As 
% \begin{array}{rcl}  C([0,1], \mathbb R) & \rightarrow &  \mathbb R\\
$\varphi  \mapsto  \varphi(1)$
%\end{array}$
is continuous, by contraction principle (\cite[Theorem 4.2.1]{DZ}), we get that $\lambda_1(1)$
satisfies  a LDP  with good rate function 
$ J_{0,\theta}, $ 
where, for any $\eta \in [0,1[,$ we denote by
 $$ I_\eta(\varphi) = \int_\eta^1 f(t, \varphi(t), \dot\varphi(t)) dt,$$
with
$$f(t, x,y) = \frac{1}{2}\left( y - \frac{1}{2t}\left(x - \sqrt{x^2 -4t}\right) \right)^2$$
and, for $x \geq 2,$ $\theta \geq 2\sqrt \eta,$
$$ J_{\eta,\theta}(x) = \inf_{\substack{
\varphi \, \textrm{s.t.} \, \varphi(\eta)=\theta,\\
 \varphi(1) = x}} I_\eta(\varphi).$$

As $I_\eta(\varphi)$ is a good rate function, the infimum in the above problem is reached,
and we denote by $\varphi^\eta$ an infimum.
 For a smooth function $f,$ the classical theory of extremal problems (see eg. \cite{Ioffe})
predicts that $\varphi^\eta$ should be solution of the Euler-Lagrange equation \eqref{Euler}.
Of course, the constraint $\varphi(t) \ge 2\sqrt t$ will play a crucial role and because the lack of smoothness
of $f$ around $x= 2\sqrt t,$ the situation is more involved. However, most of the arguments used in the sequel are 
classical in such a context and we will only give sketches of the proofs.

\smallskip
We first show that the solution of    Euler-Lagrange equation \eqref{Euler} realises the infimum
among functions staying strictly above $t \mapsto 2\sqrt t.$ Namely we have
\begin{lemma} 
\label{opt1}
For any $\eta \in ]0,1[,$ if, for any $t \in [\eta, 1],$ $\varphi^\eta(t) >2\sqrt t,$ 
then, $$\varphi^\eta(t) = \frac{x-\theta}{1- \eta}(t-\eta) + \theta.$$
\end{lemma}

\noindent
{\bf Proof:} Let $\eta \in [0,1[$
 be fixed. % We denote by $\epsilon := \inf_{t \in [\eta,1]} (\varphi^\eta(t) -2\sqrt t)>0.$\\

 It is easy to check that the infimum of $I_\eta$ is finite, therefore, we know that
$\varphi^\eta$ is absolutely continuous with $\dot \varphi^\eta \in \mathbb L^1.$
 Following the proof of Theorem 4 in Chapter 9.2.3 in \cite{Ioffe} (the details are left to the reader), one can show that
it is a solution to the DuBois-Reymond equation, i.e. there exists a constant $r$ such that
for any $t \in [\eta,1],$
\begin{equation}\label{DBR} 
 \frac{\partial f}{\partial y} (t, \varphi^\eta(t), \dot \varphi^\eta(t))- \int_\eta^t \frac{\partial f}{\partial x}(s, \varphi^\eta(s), \dot \varphi^\eta(s)) =r.
\end{equation}

%For $m \in \mathbb N^*,$ let 
%$$ \Delta_m = \{t \in [\eta,1], |\dot\varphi^\eta(t) |\leq m\}$$

 From there, one can check that $\varphi^\eta$ is a $\mathcal C^2$ solution of the Euler-Lagrange equation, that is, for any $t \in [\eta, 1],$
\begin{equation}\label{Euler}
\frac{d}{dt} \frac{\partial f}{\partial y}(t, \varphi^\eta(t), \dot \varphi^\eta(t))- \frac{\partial f}{\partial x}(t, \varphi^\eta(t), \dot \varphi^\eta(t)) =0.
\end{equation}

Indeed, if we define
$$ g(t,y) := f(t,\varphi^\eta(t),y )-y \int_\eta^t \frac{\partial f}{\partial x}(s, \varphi^\eta(s), \dot \varphi^\eta(s))ds - ry,$$
%With our expression of $f,$
 $g(t,.)$ is a convex quadratic polynomial, therefore it has a unique minimizer $y(t)$,
which is solution of the equation
$$ \frac{\partial g}{\partial y} (t,y(t)) = 0.$$
We can compute $y$ explicitely, it is given by 
$$ y(t) =  \frac{\varphi^\eta(t)}{2t} - \frac{1}{2t} \sqrt{\varphi^\eta(t)^2-4t}
+ \frac{1}{2t}\frac{\varphi^\eta(t)}{\sqrt{\varphi^\eta(t)^2-4t}} +r.$$
As we know that $\varphi^\eta$ is absolutely continuous, so is $y.$
But, by unicity of the minimizer, we have that $y = \dot \varphi^\eta$ so that we get
that $\varphi^\eta$ is continuously differentiable. Therefore $\frac{\partial g}{\partial y}$
is continuously differentiable in both variables and from the implicit function theorem,
we get that $y(t)$ is continuously differentiable, so that $\varphi^\eta$ is twice continuously differentiable.
 Differentiating \eqref{DBR} we get \eqref{Euler}.

%The proof follows the same lines as in the second part of the proof of 
%Theorem 4 in Chapter 9.2.3 in \cite{Ioffe}, the only point to check is that
%the exists an $L^1$ 

Straightforward computations  leads to $\ddot{\varphi^\eta} \equiv 0$ and thus
$$\varphi^\eta(t) = \frac{x-\theta}{1- \eta}(t-\eta) + \theta. $$
\hfill $\Box$\smallskip

Now, the goal will be to determine in which cases this linear solution realises
the infimum over all admissible functions and in which cases we can do better by touching
the wall $t \mapsto 2 \sqrt t.$\smallskip

We need the comparison
\begin{lemma}\label{opt2}
 For any $\eta \in [0,1[,$ if there exists $t_0 \in [\eta,1]$
such that $\varphi^\eta(t_0) = 2\sqrt{t_0},$
then $$I_\eta(\varphi^\eta) \geq \int_2^x \sqrt{u^2-4} \ du.$$
\end{lemma}
 
\noindent
{\bf Proof:} We have
$$I_\eta(\varphi^\eta) \geq \int_{t_0}^1 f(t, \varphi^\eta(t),\dot \varphi^\eta(t))dt.$$
 For any $\varphi$ such that $I_\eta(\varphi) <\infty,$ for $t\geq \eta,$ we denote
$$ K_t(\varphi) = \frac{1}{2} \int_t^1 \left( \dot \varphi(s) - \frac{1}{2s}
(\varphi(s) + \sqrt{\varphi^2(s)-4s})\right)^2 ds.$$
If we let $y(s) = \frac{\varphi(s)}{\sqrt{s}},$ one has
$$ K_t(\varphi) = I_t(\varphi) - \int_{y(t)}^{y(1)} \sqrt{u^2-4} \,\,du.$$
But $K_t \geq 0$ so that
$$I_\eta(\varphi^\eta) \geq \int_{\frac{\varphi^\eta(t_0)}{\sqrt{t_0}}}^{x} \sqrt{u^2-4} \,\,du
= \int_2^x \sqrt{u^2-4} du.$$
\hfill $\Box$\smallskip

 From there, we can prove 
\begin{lemma}\label{opt3}
 For any $\eta \in [0,1[,$ 

\noindent if $\theta = 2\sqrt \eta$ or  $\left(2\sqrt \eta < \theta < 1+\eta
\textrm{ and } x \leq \frac{\theta + \sqrt{\theta^2-4\eta}}{2}+ \frac{2}{\theta + \sqrt{\theta^2-4\eta}}\right),$
then $$ J_{\eta, \theta}(x) = \int_2^x \sqrt{u^2-4} \,\,du.$$
\end{lemma}

\noindent
{\bf Proof:} The proof consists in exhibiting explicit functions $\varphi_\eta^*$
realizing the infimum, that is such that $ I_{\eta}(\varphi_\eta^*) = \int_2^x \sqrt{u^2-4} du.$

\noindent For $\theta = 2\sqrt \eta$ and $x \geq 2,$ we let $t^*$ be such that  $\sqrt{t^*}:= \frac{x+\sqrt{x^2-4}}{2}$ and
\[\varphi^*_\eta(t) = \left\{ 
\begin{array}{ll}
 2\sqrt t & \textrm{ if } \eta \leq t \leq t^* \\
2\sqrt{t^*} + \frac{1}{\sqrt{t^*}}(t-t^*) & \textrm{ if } t^*\vee \eta \leq t\leq 1
\end{array}
\right.\]
 For $2\sqrt \eta < \theta < 1+\eta$
and $x \leq \frac{\theta + \sqrt{\theta^2-4\eta}}{2}+ \frac{2}{\theta + \sqrt{\theta^2-4\eta}},$
we let $s^*$ be such that $\sqrt{s^*} = \frac{\theta+\sqrt{\theta^2-4\eta}}{2},$ and
\[\varphi^*_\eta(t) = \left\{ 
\begin{array}{ll}
\theta + \frac{1}{\sqrt{s^*}}(t-\eta) & \textrm{ if } \eta \leq t \leq s^* \\
 2\sqrt t & \textrm{ if } s^* \leq t \leq t^* \\
2\sqrt{t^*} + \frac{1}{\sqrt{t^*}}(t-t^*) & \textrm{ if } t^* \leq t\leq 1
\end{array}
\right.\]
The extension to the case when $\eta=0$ is easy to obtain and left to the reader.
\hfill $\Box$ \smallskip

In the other cases, tedious computations allows to check that the solution of Euler-Lagrange
realises the infimum. More precisely we have 
\begin{lemma} \label{opt4}
 For any $\eta \in ]0,1[,$ 

\noindent if $\theta > 1+\eta$ or $\left(2\sqrt \eta < \theta < 1+\eta
\textrm{ and } x >\frac{\theta + \sqrt{\theta^2-4\eta}}{2}+ \frac{2}{\theta + \sqrt{\theta^2-4\eta}}\right),$
then $$ J_{\eta, \theta}(x) = I_\eta(d_{\eta,\theta,x}),$$
where for any $t \in [\eta, 1],$
$$ d_{\eta,\theta,x}(t) = \frac{(x-\theta)}{1-\eta}(t-\eta) +\theta.$$
\end{lemma}

The last point, to complete the proof of the Theorem is to extend the above lemma to the case when $\eta = 0.$ More precisely, we have,
\begin{lemma}\label{opt5}
 If $\theta \geq 1$ or $\left(0<\theta <1 \textrm{  and } x \geq \theta + \frac{1}{\theta}\right),$
let $\varphi_*(t) = (x-\theta)t + \theta$ for $t \in [0,1],$ then for any 
$\varphi,$  $I_0(\varphi_*) \leq I_0(\varphi)$
\end{lemma}
We do not detail the proof of the lemma which consists in cutting the integral defining $I_\theta,$
on $[0,\eta]$ with $\eta $ small enough for the integral to be small and on $[\eta, 1]$
on which the previous lemmas give the minimizers.\smallskip

This  concludes the proof by the following computation.
We set $\alpha = x - \theta$ and use the change of variable  $u = \theta + \alpha t + \sqrt{(\theta + \alpha t )^2 -4t}.$
\begin{eqnarray*}
I(\varphi_*) &=& \frac{1}{2} \int_0^1 \left( \alpha - \frac{1}{2t}( \theta + \alpha t - \sqrt{(\theta + \alpha t )^2 -4t}) \right)^2 dt \\
&=&  - \ln\left( \frac{x + \sqrt{x^2 -4}}{2 \theta} \right)  + \frac{1}{4} x \sqrt{x^2 -4} + \frac{1}{4} x^2 - \theta x + \frac{\theta^2}{2} + \frac{1}{2}.
 \end{eqnarray*}
 This agrees with the formulae giving $M_\theta$ and $L_\theta$ since we have
$$  \int_2^x \sqrt{z^2 - 4}\, dz = - 2\ln\left( \frac{x + \sqrt{x^2 -4}}{2 \theta} \right)  + \frac{1}{2} x \sqrt{x^2 -4}.  \quad \quad \quad \quad \quad \quad \Box$$\smallskip

\section*{Acknowledgements}
We would like to thank Marc Yor who suggested this problem to us. We are really indebted to  Ofer Zeitouni, who helped us with a crucial suggestion leading to the argument for the upper bound developed in Section 5.  We thank for his hospitality the American Institute of Mathematics (AIM) where the discussion with OZ took place.

\bibliographystyle{alpha}
\bibliography{09-21}

\end{document}